\documentclass{amsart}
\usepackage{times,mathrsfs,multicol,multibib}
\newcites{AM,MS}{References,References}
\usepackage{cite}

\title{Seminormal forms and Gram determinants for cellular algebras}
\author{Andrew Mathas}

\swapnumbers
\numberwithin{equation}{section}
\newtheorem{Defn}[equation]{Definition}
\newtheorem{Theorem}[equation]{Theorem}
\newtheorem{Prop}[equation]{Proposition}
\newtheorem{Lemma}[equation]{Lemma}
\newtheorem{Cor}[equation]{Corollary}

\theoremstyle{remark}
\newtheorem{Point}[equation]{}
\newtheorem{Remark}[equation]{Remark}

\newenvironment{Example}[1][\relax]%
  {\refstepcounter{equation}
   \trivlist%
   \ifx#1\relax
     \item[\hskip\labelsep\theequation\space\textbf{Example}]
   \else
     \item[\hskip\labelsep\theequation\space\textbf{Example} (#1)]
   \fi
   \ignorespaces
  }{\unskip\nobreak\hfil%
    \penalty50\hskip2em\hbox{}\nobreak\hfil$\Diamond$%
    \parfillskip=0pt\finalhyphendemerits=0\penalty-100\endtrivlist}

\def\B{\mathscr{B}_n(q,r)}
\def\CC{\mathscr C}
\def\RR{\mathscr R}

\def\H{\mathscr H}
\def\L{\mathscr L}

\def\Sym{\mathfrak S}
\def\T{\mathbb T}
\def\Z{\mathbb Z}

\def\s{\mathsf s}
\def\t{\mathsf t}
\def\fS{\mathsf S}
\def\fT{\mathsf T}
\def\u{\mathsf u}
\def\v{\mathsf v}
\def\w{\mathsf w}
\def\x{\mathsf x}
\def\y{\mathsf y}
\def\z{\mathsf z}

\let\bar=\overline
\let\gdom=\vartriangleright
\let\gedom=\trianglerighteq
\def\map#1#2{\,{:}\,#1\!\longrightarrow\!#2}
\def\floor#1#2{\lfloor\tfrac{#1}{#2}\rfloor}

\let\<=\langle
\let\>=\rangle
\def\({\Big(}
\def\){\Big)}
\def\Prod{\displaystyle\prod}

\def\And{\text{ and } }
\def\For{\text{ for } }
\def\ForAll{\text{ for all } }
\def\If{\text{ if } }
\def\Otherwise{\text{ otherwise} }

\DeclareMathOperator{\End}{End}
\DeclareMathOperator{\rad}{rad}

\DeclareMathOperator{\Span}{Span}

{\catcode`\|=\active
  \gdef\set#1{\mathinner{\lbrace\,{\mathcode`\|"8000%
                                   \let|\midvert #1}\,\rbrace}}
}
\def\midvert{\egroup\mid\bgroup}
{\catcode`\|=\active
  \gdef\Set#1{\mathinner{\Big\lbrace\,{\mathcode`\|"8000%
                                   \let|\Midvert #1}\,\Big\rbrace}}
}
\def\Midvert{\egroup\,\Big|\,\bgroup}

\newcommand{\thi}[1]{\mbox{\rm \underline{\textbf{#1}}}}
\newcommand{\mis}[1]{\mbox{$\hat{{\rm\underline{\textbf{#1}}}}$}}
\newcommand{\ul}[1]{\underline{#1}}

\DeclareMathOperator{\ecal}{\mathcal{E}}
\DeclareMathOperator{\hcal}{\mathcal{H}}
\DeclareMathOperator{\lcal}{\mathcal{L}}
\DeclareMathOperator{\ncal}{\mathcal{N}}
\DeclareMathOperator{\scal}{\mathcal{S}}
\DeclareMathOperator{\tcal}{\mathcal{T}}
\DeclareMathOperator{\mfrk}{\mathfrak{m}}
\DeclareMathOperator{\ufrk}{\mathfrak{u}}
\DeclareMathOperator{\ebb}{\mathbb{E}}
\DeclareMathOperator{\fbb}{\mathbb{F}}
\DeclareMathOperator{\lbb}{\mathbb{L}}
\DeclareMathOperator{\zbb}{\mathbb{Z}}

\DeclareMathOperator{\Mat}{Mat}
\DeclareMathOperator{\supp}{supp}
\DeclareMathOperator{\One}{\textbf{1}}

\begin{document}
\makeatletter
\begin{center}\baselineskip14pt\relax
{\bfseries\uppercasenonmath\@title\@title}\\
\bigskip
\textsc{Andrew Mathas}\\
\bigskip
\textit{School of Mathematics and Statistics F07, 
        University of Sydney,\\ Sydney NSW 2006, Australia}\\
E-mail address: \texttt{a.mathas@usyd.edu.au}
\end{center}
\markboth{\uppercase{Andrew Mathas}}%
         {\uppercase{Seminormal forms and Gram determinants}}

\bibliographystyleAM{andrew}

\begin{abstract}
This paper develops an abstract framework for constructing
``seminormal forms'' for cellular algebras. That is, given a cellular
$R$--algebra $A$ which is equipped with a \textit{family of
JM--elements} we give a general technique for constructing orthogonal
bases for~$A$, and for all of its irreducible representations, when
the JM--elements \textit{separate}~$A$. The seminormal forms for $A$
are defined over the field of fractions of~$R$.  Significantly, we
show that the Gram determinant of each irreducible $A$--module is
equal to a product of certain structure constants coming from the
seminormal basis of~$A$. In the non--separated case we use our
seminormal forms to give an explicit basis for a block decomposition
of~$A$. 
\end{abstract}

\@setabstract
\makeatother

\section{Introduction}

The purpose of this paper is to give an axiomatic way to construct
``seminormal forms'' and to compute Gram determinants for the
irreducible representations of semisimple cellular algebras. By this
we mean that, starting from a given cellular basis
$\{a^\lambda_{\s\t}\}$ for a cellular algebra $A$, we give a new
cellular basis $\{f^\lambda_{\s\t}\}$ for the algebra which is
orthogonal with respect to a natural bilinear form on the algebra.
This construction also gives a ``seminormal basis'' for each of the
cell modules of the algebra. We show that the Gram determinant of the
cell modules (the irreducible $A$--modules) can be computed in terms
of the structure constants of the new cellular basis of $A$. Combining
these results gives a recipe for computing the Gram determinants 
of the irreducible $A$--modules.

Of course, we cannot carry out this construction for an arbitrary
cellular algebra~$A$. Rather, we assume that the cellular algebra
comes equipped with a family of ``Jucys--Murphy'' elements. These are
elements of $A$ which act on the cellular basis of~$A$ via upper
triangular matrices. We will see that, over a field, the existence of
such a basis $\{f^\lambda_{\s\t}\}$ forces $A$ to be (split)
semisimple. The cellular algebras which have JM--elements include the
group algebras of the symmetric groups, any split semisimple algebra,
the Hecke algebras of type~$A$, the $q$--Schur algebras, the
(degenerate) Ariki--Koike algebras, the cyclotomic $q$--Schur
Algebras, the Brauer algebras and the BMW algebras.

At first sight, our construction appears to be useful only in the
semisimple case. However, in the last section of this paper we apply
these ideas in the non--semisimple case to construct a third cellular
basis $\{g^\lambda_{\s\t}\}$ of~$A$. We show that this basis gives an
explicit decomposition of $A$ into a direct sum of smaller cellular
subalgebras.  In general, these subalgebras need not be
indecomposable, however, it turns out that these subalgebras are
indecomposable in many of the cases we know about. As an application,
we give explicit bases for the block decomposition of the group
algebras of the symmetric groups, the Hecke algebras of type~$A$, the
Ariki--Koike algebras with $q\ne1$, the degenerate Ariki--Koike
algebras and the (cyclotomic) $q$--Schur algebras.

There are many other accounts of seminormal forms in the literature;
see, for example, \citeAM{ABR:I,Hoefsmit,LeducRam,Ram:seminormal}. The
main difference between this paper and previous work is that, starting
from a cellular basis for an algebra we \textit{construct} seminormal
forms for the entire algebra, rather than just the irreducible
modules. The main new results that we obtain are explicit formulae for
the Gram determinants of the cell modules in the separated case, and a
basis for a block decomposition of the algebra in the non--separated
case. These seminormal forms that we construct have the advantage that
they are automatically defined over the field of fractions of the base
ring; this is new for the Brauer and BMW algebras.

It follows using the Wedderburn theorem that an algebra has a family
of separating JM--elements if and only if it is split semisimple (see
Example~\ref{semisimple ex}). As every split semisimple algebra is
cellular this suggests that cellular algebras provide the right
framework for studying seminormal forms.  There is, however, an
important caveat: the set of JM--elements for a cellular algebra is not
canonical as it  depends heavily on the particular choice of cellular
basis. Consequently, to study an algebra using the techniques in this
paper one has to first find a cellular basis for the algebra and then
find an appropriate set of JM--elements. Neither of these tasks is
necessarily easy especially as, ideally, we would like the
set of JM--elements to be compatible with modular reduction.

In the appendix to this paper, Marcos Soriano, gives an alternative
matrix theoretic approach to the theory of seminormal forms. Using
only the Cayley--Hamilton theorem he shows that if you have a family
of operators acting on a module via upper triangular matrices which
satisfy an analogous separation condition then you can construct a
complete set of pairwise orthogonal idempotents. This shows that,
ultimately, the theory of seminormal forms rests on the
Cayley--Hamilton theorem. Note that unlike in our treatment, Soriano
does not need to assume that the JM--elements commute or that they are
$*$--invariant.

This paper is organized as follows. In the next section we recall
Graham and Lehrer's theory of cellular algebras and define
JM--elements for cellular algebras. We then show that any cellular
algebra with a family of separating JM--elements is necessarily
semisimple and, by way of example, show that most of the well--known
cellular algebras have JM--elements. The third section of the paper
develops the theory of JM--elements in the separated case, culminating
with the construction of a seminormal basis for a cellular algebra and
the computation of the Gram determinants of the cell modules. In the
last section of the paper we use modular reduction to study the
non--separated case. Our main result gives a cellular basis for a
decomposition of the original cellular algebra into blocks.  Finally,
in the appendix Marcos Soriano gives his matrix theoretic approach
to the theory of seminormal forms.

\section{Cellular algebras and JM--elements}

We begin by recalling Graham and Lehrer's~\citeAM{GL} definition of a
cellular algebra. Let $R$ be commutative ring with~$1$ and let $A$ be
a unital $R$--algebra and let $K$ be the field of fractions of $R$.

\begin{Defn}[Graham and Lehrer]\label{cellular}
A \textsf{cell datum} for $A$ is a triple $(\Lambda,T,C)$
where $\Lambda=(\Lambda,>)$ is a finite poset,
$T(\lambda)$ is a finite set for each $\lambda\in\Lambda$, and 
$$C\map{\coprod_{\lambda\in\Lambda}T(\lambda)\times T(\lambda)}
                A; (\s,\t)\mapsto a^\lambda_{\s\t}$$
is an injective map (of sets) such that:
\begin{enumerate}
\item $\set{a^\lambda_{\s\t}|\lambda\in\Lambda, \s,\t\in T(\lambda)}$
is an $R$--free basis of $A$;
\item For any $x\in A$ and $\t\in T(\lambda)$ there exist
scalars $r_{\t\v x}\in R$ such that, for any $\s\in T(\lambda)$,
$$a^\lambda_{\s\t}x\equiv 
\sum_{\v\in T(\lambda)} r_{\t\v x}a^\lambda_{\s\v} \pmod{A^\lambda},$$
where $A^\lambda$ is the $R$--submodule of $A$ with basis
$\set{a^\mu_{\y\z}|\mu>\lambda\And\y,\z\in T(\mu)}$.
\item The $R$--linear map determined by $*\map AA;
a^\lambda_{\s\t}=a^\lambda_{\t\s}$, for all $\lambda\in\Lambda$ and
$\s,\t\in T(\lambda)$, is an anti--isomorphism of $A$.
\end{enumerate}
If a cell datum exists for $A$ then we say that $A$ is a
\textsf{cellular algebra}. 
\end{Defn}

Henceforth, we fix a cellular algebra $A$ with cell datum
$(\Lambda,T,C)$ as above. We will also assume that $T(\lambda)$ is a
poset with ordering $\gdom_\lambda$, for each $\lambda\in\Lambda$.
For convenience we set
$T(\Lambda)=\coprod_{\lambda\in\Lambda}T(\lambda)$. We consider
$T(\Lambda)$ as a poset with the ordering $\s\gdom\t$ if either (1)
$\s,\t\in T(\lambda)$, for some $\lambda\in\Lambda$, and
$\s\gdom_\lambda\t$, or (2) $\s\in T(\lambda)$, $\t\in T(\mu)$ and
$\lambda>\mu$.  We write $\s\gedom\t$ if $\s=\t$ or $\s\gdom\t$. If
$\s\gedom\t$ we say that $\s$ \textsf{dominates} $\t$.

Note that, by assumption $A$, is a free $R$--module of finite rank
$|T(\Lambda)|$.

Let $A_K=A\otimes_R K$. As $A$ is free as an $R$--module, $A_K$ is a
cellular algebra with cellular basis 
$\set{a^\lambda_{\s\t}\otimes 1_K|\lambda\in\Lambda\And\s,\t\in T(\lambda)}$. 
We consider $A$ as a subalgebra of $A_K$ and, abusing notation,
we also consider $a^\lambda_{\s\t}$ to be elements of~$A_K$.

We recall some of the general theory of cellular
algebras. First, applying the $*$~involution to part~(b) of
Definition~\ref{cellular} we see that
if $y\in A$ and $\s\in T(\lambda)$ then there exist scalars 
$r_{\s\u y}'\in R$ such that, for all $\t\in T(\lambda)$,
\begin{equation}\label{left mult}
ya^\lambda_{\s\t}\equiv 
     \sum_{\u\in T(\lambda)} r_{\s\u y}'a^\lambda_{\u\t} \pmod {A^\lambda}.
\end{equation}
Consequently, $A^\lambda$ is a two--sided ideal of $A$, for any
$\lambda\in\Lambda$.

Next, for each $\lambda\in\Lambda$ define the \textsf{cell module}
$C(\lambda)$ to be the free $R$--module with basis
$\set{a^\lambda_\t|\t\in T(\lambda)}$ and with
$A$--action given by 
$$a^\lambda_\t x= 
       \sum_{\v\in T(\lambda)} r_{\t\v x}a^\lambda_{\v},$$
where $r_{\t\v x}$ is the same scalar which appears in
Definition~\ref{cellular}. As $r_{\t\v x}$ is independent of~$\s$ this
gives a well--defined $A$--module structure on $C(\lambda)$. The map
$\<\ ,\ \>_\lambda\map{C(\lambda)\times C(\lambda)}R$
which is determined by
\begin{equation}\label{inner product}
\<a^\lambda_\t,a^\lambda_\u\>_\lambda a^\lambda_{\s\v}
     \equiv a^\lambda_{\s\t}a^\lambda_{\u\v}\pmod{A^\lambda},
\end{equation}
for $\s, \t, \u, \v\in T(\lambda)$, defines a symmetric bilinear form 
on $C(\lambda)$. This form is associative in the sense that
$\<ax,b\>_\lambda=\<a,bx^*\>_\lambda$, for all $a,b\in C(\lambda)$ and
all $x\in A$. From the definitions, for any $\s\in T(\lambda)$ the
cell module $C(\lambda)$ is naturally isomorphic to the $A$--module
spanned by $\set{a^\lambda_{\s\t}+A^\lambda|\t\in T(\lambda)}$. The
isomorphism is the obvious one which sends 
$a^\lambda_\t\mapsto a^\lambda_{\s\t}+A^\lambda$, for 
$\t\in T(\lambda)$.

For $\lambda\in\Lambda$ we define 
$\rad C(\lambda)=\set{x\in C(\lambda)|\<x,y\>_\lambda=0
         \ForAll y\in C(\lambda)}$. As the bilinear form on
$C(\lambda)$ is associative it follows that $\rad C(\lambda)$ is an
$A$--submodule of $C(\lambda)$. Graham and
Lehrer~\citeAM[Theorem~3.4]{GL} show that the $A_K$--module
$D(\lambda)=C(\lambda)/\rad C(\lambda)$ is absolutely irreducible and
that $\set{D(\lambda)\ne0|\lambda\in\Lambda}$ is a
complete set of pairwise non--isomorphic irreducible $A_K$--modules.

The proofs of all of these results follow easily from
Definition~\ref{cellular}. For the full details see
\citeAM[\S2--3]{GL} or \citeAM[Chapt.~2]{M:Ulect}.

In this paper we are interested only in those cellular algebras which
come equipped with the following elements.

\begin{Defn}\label{JM}
A family of \textsf{JM--elements} for $A$ is a set $\{L_1,\dots,L_M\}$
of commuting elements of $A$ together with a set of scalars,
$\set{c_\t(i)\in R|\t\in T(\Lambda)\And1\le i\le M}$, such that 
for $i=1,\dots,M$ we have $L_i^*=L_i$ and, for
all $\lambda\in\Lambda$ and $\s,\t\in T(\lambda)$,
$$ a^\lambda_{\s\t}L_i
   \equiv c_\t(i)a^\lambda_{\s\t}
        +\sum_{\v\gdom\t}r_{\t\v}a^\lambda_{\s\v} \pmod{A^\lambda},
$$
for some $r_{\t\v}\in R$ $($which depend on~$i)$. 
We call $c_\t(i)$ the \textsf{content} of $\t$ at $i$.
\end{Defn}

Implicitly, the JM--elements depend on the choice of cellular basis
for~$A$.

Notice that we also have the following left hand
analogue of the formula in (\ref{JM}):
\begin{equation}\label{JM*}
L_ia^\lambda_{\s\t}
   \equiv c_\s(i)a^\lambda_{\s\t}
        +\sum_{\u\gdom\s}r_{\s\u}'a^\lambda_{\u\t} \pmod{A^\lambda},
\end{equation}
for some $r_{\s\u}'\in R$. 

\begin{Point}\label{L_K}
Let $\L_K$ be the subalgebra of $A_K$ which is generated by
$\{L_1,\dots,L_M\}$. By definition,~$\L_K$ is a commutative subalgebra
of~$A_K$. It is easy to see that each $\t\in T(\Lambda)$ gives rise to
a one dimensional representation $K_\t$ of~$\L_K$ on which $L_i$ acts
as multiplication by $c_\t(i)$, for $1\le i\le M$.  In fact, since
$\L_K$ is a subalgebra of $A_K$, and $A_K$ has a filtration by cell
modules, it follows that $\set{K_\t|\t\in T(\Lambda)}$ is a complete
set of irreducible $\L_K$--modules.
\end{Point}

These observations give a way of detecting when $D(\lambda)\ne0$
(cf.~\citeAM[Prop.~5.9(i)]{GL}). 

\begin{Prop}\label{simple crit}
Let $A_K$ be a cellular algebra with a family of JM--elements and fix
$\lambda\in\Lambda$, and $\s\in T(\lambda)$. Suppose that whenever
$\t\in T(\Lambda)$ and $\s\gdom\t$ then $c_\t(i)\ne c_\s(i)$, for 
some~$i$ with $1\le i\le M$. Then $D(\lambda)\ne0$.
\end{Prop}

\begin{proof}
By definition~\ref{JM}, for any $\mu\in\Lambda$ the $\L_K$--module
    composition factors of $C(\mu)$ are precisely the modules
    $\set{K_\s|\s\in T(\mu)}$. Observe that if $\u,\v\in T(\Lambda)$
    then $K_\u\cong K_\v$ as $\L_K$--modules if and only if 
    $c_\u(i)=c_\v(i)$, for $1\le i\le M$. Therefore, our assumptions 
    imply that $K_\t$ is not an $\L_K$--module composition factor of any
    cell module $C(\mu)$ whenever $\lambda>\mu$. Consequently, $K_\t$ is not
    an $\L_K$--module composition factor of $D(\mu)$ whenever
    $\lambda>\mu$. However, by \citeAM[Prop.~3.6]{GL}, $D(\mu)$ is
    a composition factor of $C(\lambda)$ only if $\lambda\ge\mu$.
    Therefore, $a^\lambda_\t\notin\rad C(\lambda)$ and,
    consequently, $D(\lambda)\ne0$ as claimed.
\end{proof}

Motivated by  Proposition~\ref{simple crit}, we break our study of
cellular algebras with JM--elements into two cases depending upon
whether or not the condition in Proposition~\ref{simple crit} is
satisfied.

\begin{Defn}[Separation condition]\label{separation}
Suppose that $A$ is a cellular algebra with JM--elements
$\{L_1,\dots,L_M\}$. The JM--elements \textsf{separate} $T(\lambda)$
$($over $R)$ if whenever $\s,\t\in T(\Lambda)$ and $\s\gdom\t$ then
$c_\s(i)\ne c_\t(i)$, for some $i$ with $1\le i\le M$.
\end{Defn}

In essence, the separation condition says that the contents $c_\t(i)$
distinguish between the elements of $T(\Lambda)$. Using the argument of
Proposition~\ref{simple crit} we see that the separation condition
forces $A_K$ to be semisimple.

\begin{Cor}\label{semisimple}
Suppose that $A_K$ is a cellular algebra with a family of JM--elements
which separate $T(\Lambda)$. Then $A_K$ is (split) semisimple.
\end{Cor}

\begin{proof}
By the general theory of cellular algebras~\citeAM[Theorem~3.8]{GL},
$A_K$ is (split) semisimple if and only if $C(\lambda)=D(\lambda)$ for
all $\lambda\in\Lambda$. By the argument of 
Proposition~\ref{simple crit}, the separation condition implies that
if $\t\in T(\lambda)$ then $K_\t$ does not occur as an $\L_K$--module
composition factor of $D(\mu)$ for any $\mu>\lambda$. By
\citeAM[Prop.~3.6]{GL}, $D(\mu)$ is a composition factor of
$C(\lambda)$ only if $\lambda\ge\mu$, so the cell module
$C(\lambda)=D(\lambda)$ is irreducible. Hence, $A_K$ is semisimple as
claimed.
\end{proof}

In Example~\ref{semisimple ex} below we show that every split
semisimple algebra is a cellular algebra with a family of JM--elements
which separate~$T(\Lambda)$.

\begin{Remark}
Corollary~\ref{semisimple} says that if a cellular algebra $A$ has a
family of JM--elements which separate $T(\Lambda)$ then $A_K$ is split
semisimple.  Conversely, we show in Example~\ref{semisimple ex} below
that every split semisimple algebra has a family of JM--elements which
separate $T(\Lambda)$. However, if $A$ is semisimple and $A$ has a
family of JM--elements then it is not true that the JM--elements must
separate~$A$; the problem is that an algebra can have different
families of JM--elements. As described in Example~\ref{BMW ex} below,
the Brauer and BMW algebras both have families of JM--elements. Combined
with work of Enyang \citeAM[Examples~7.1 and~10.1]{Enyang:ss} this
shows that there exist BMW and Brauer algebras which are semisimple
and have JM--elements which do not separate~$T(\Lambda)$.
\end{Remark}

\begin{Remark} Following ideas of
Grojnowski~\citeAM[(11.9)]{Klesh:book} and (\ref{L_K}) we can use the
algebra $\L_K$ to define \textit{formal characters}
of $A_K$--modules as follows. Let $\set{K_\t|\t\in L(\Lambda)}$ be a
complete set of non--isomorphic irreducible $\L_K$--modules, where
$L(\Lambda)\subseteq T(\Lambda)$. If~$M$ is any $A_K$--module let
$[M:K_\t]$ be the decomposition multiplicity of the irreducible
$\L_K$--module $K_\t$ in $M$. Define the formal character of~$M$ to be
$$\operatorname{ch}M=\sum_{\t\in L(\Lambda)}[M:K_\t]\,e^\t,$$ 
which is element of the free $\Z$--module with basis 
$\set{e^\t|\t\in L(\Lambda)}$. It would be interesting to know to what
extent these characters determine the representations of~$A$.
\end{Remark}
    
We close this introductory section by giving examples of cellular
$R$--algebras which have a family of JM--elements. Rather than
starting with the simplest example we start with the motivating
example of the symmetric group. The latter examples are either less
well--known or new.

\begin{Example}[Symmetric groups]
The first example of a family of JM--elements was given
by Jucys~\citeAM{Jucys} and, independently, by
Murphy~\citeAM{murphy:ops}. (In fact, these elements first appear in
the work of Young~\citeAM{QSAII}.) Let $A=R\Sym_n$ be the group ring
of the symmetric group of degree~$n$. Define
$$L_i=(1,i)+(2,i)+\dots+(i-1,i),\qquad\For i=2,\dots,n.$$
Murphy~\citeAM{murphy:ops} showed that these elements commute and he
studied the action of these elements on the seminormal basis of the
Specht modules. The seminormal basis of the Specht modules can be
extended to a seminormal basis of $R\Sym_n$, so Murphy's work shows
that the group algebra of the symmetric group fits into our general
framework. We do not give further details because a better approach to
the symmetric groups in given by the special case $q=1$ of
Example~\ref{HeckeA ex} below which concerns the Hecke algebra of
type~$A$. 
\end{Example}

\begin{Example}[Semisimple algebras]\label{semisimple ex}
By Corollary~\ref{semisimple} every cellular algebra over a field
which has a family of JM--elements which separate $T(\Lambda)$ is
split semisimple. In fact, the converse is also true. Note that
a cellular algebra is semisimple if and only if it is split
semisimple, so non--split semisimple algebras do not arise in our
setting. In fact, the appendix shows that in the separated case the
existence of family of JM--elements acting on a module forces absolute
irreducibilty, so JM--elements never arise in the non--split case.

Suppose that~$A_K$ is a split semisimple algebra. Then the Wedderburn
basis of matrix units in the simple components of~$A_K$ is a cellular
basis of~$A_K$.  We claim that~$A_K$ has a family of JM--elements.  To
see this it is enough to consider the case when $A_K=\text{Mat}_n(K)$
is the algebra of $n\times n$ matrices over $K$. Let $e_{ij}$ be the
elementary matrix with has a~$1$ in row~$i$ and column~$j$ and zeros
elsewhere. Then it is easy to check that $\{e_{ij}\}$ is a cellular
basis for $A_K$ (with $\Lambda=\{1\}$, say, and
$T(\lambda)=\{1,\dots,n\})$). Let $L_i=e_{ii}$, for $1\le i\le n$.
Then $\{L_1,\dots,L_n\}$ is a family of JM--elements for~$A_K$ which
separate $T(\Lambda)$.

By the last paragraph, any split semisimple algebra $A_K$ has a family
of JM--elements $\{L_1,\dots,L_M\}$ which separate~$T(\Lambda)$, where
$M=d_1+\dots+d_r$ and $d_1,\dots,d_r$ are the dimensions of the
irreducible $A_K$--modules.  The examples below show that we can often
find a much smaller set of JM--elements. In particular, this shows
that the number~$M$ of JM--elements for an algebra is not an invariant
of~$A$!  Nevertheless, in the separated case we will show that the
JM--elements are always linear combinations of the diagonal elementary
matrices coming from the different Wedderburn components of the
algebra.  Further, the subalgebra of~$A_K$ generated by a family of
JM--elements is a maximal abelian subalgebra of~$A_K$.
\end{Example}

If $A_K$ is a cellular algebra and explicit formulae for the
Wedderburn basis of~$A_K$ are known then we do not need this paper to
understand the representations of~$A_K$.  One of the points of this
paper is that \textit{if we have a cellular basis for an $R$--algebra
$A$ together with a family of JM--elements then we can construct a
Wedderburn basis for~$A_K$.}

\begin{Example}[A toy example]\label{toy ex}
Let $A=R[X]/(X-c_1)\dots(X-c_n)$, where $X$ is an indeterminate over
$R$ and $c_1,\dots,c_n\in R$. Let $x$ be the image of~$X$ in~$A$ under
the canonical projection $R[X]\longrightarrow A$. Set
$a_i:=a^i_{ii}=\prod_{j=1}^{i-1}(x-c_j)$, for $i=1,\dots,n+1$. Then
$A$ is a cellular algebra with $\Lambda=\{1,\dots,n\}$, $T(i)=\{i\}$,
for $1\le i\le n$, and with cellular basis
$\{a^1_{11},\dots,a^n_{nn}\}$. Further, $x$ is a JM--element for $A$
because
$$ a_ix = (x-c_1)\dots(x-c_{i-1})x = c_ia_i+a_{i+1},$$
for $i=1,\dots,n$. Thus, $c_i(x)=c_i$, for
all~$i$. The `family' of JM--elements $\{x\}$ separates $T(\Lambda)$
if and only if $c_1,\dots,c_n$ are pairwise distinct.
\end{Example}

\begin{Example}[Hecke algebras of type~$A$]\label{HeckeA ex}
Fix an integer $n>1$ and an invertible element $q\in R$.  Let
$\H=\H_{R,q}(\Sym_n)$ be the \textsf{Hecke algebra of type $A$}. In
particular, if $q=1$ then $\H_{R,q}(\Sym_n)\cong R\Sym_n$. In general,
$\H$ is free as an $R$--module with basis $\set{T_w|w\in\Sym_n}$ and
with multiplication determined by
$$ T_{(i,i+1)}T_w=\begin{cases}
          T_{(i,i+1)w},&\If i^w>(i+1)^w,\\
         qT_{(i,i+1)w}+(q-1)T_w,&\Otherwise.
\end{cases}$$

Recall that a partition of~$n$ is a weakly decreasing sequence of
positive integers which sum to~$n$. Let $\Lambda$ be the set of
partitions of~$n$ ordered by dominance~\citeAM[3.5]{M:Ulect}. If
$\lambda=(\lambda_1,\dots,\lambda_k)$ is a partition let
$[\lambda]=\set{(r,c)|1\le c\le\lambda_r,r\le k}$ be the diagram
of~$\lambda$. A \textsf{standard $\lambda$--tableau} is a map
$\t\map{[\lambda]}\{1,\dots,n\}$ such that $\t$ is monotonic
increasing in both coordinates (i.e. rows and columns).

Given $\lambda\in\Lambda$ let $T(\lambda)$ be the set of standard
$\lambda$--tableau, ordered by dominance (the Bruhat order;
see~\citeAM[Theorem~3.8]{M:Ulect}). Murphy~\citeAM{murphy:basis} has
shown that $\H$ has a cellular basis of the form
$\set{m^\lambda_{\s\t}|\lambda\in\Lambda\And \s,\t\in T(\lambda)}$.

Set $L_1=0$ and define 
$$L_i = \sum_{j=1}^{i-1}q^{j-i}T_{(i,j)},\qquad\For 2\le i\le n.$$ 
It is a straightforward, albeit tedious, exercise to check that these
elements commute; see, for example, \citeAM[Prop.~3.26]{M:Ulect}. The
cellular algebra $*$~involution of~$\H$ is the linear extension of
the map which sends $T_w$ to $T_{w^{-1}}$, for $w\in\Sym_n$. So
$L_i^*=L_i$, for all~$i$.

For any integer $k$ let $[k]_q=1+q+\dots+q^{k-1}$ if $k\ge0$ and set
$[k]_q=-q^{-k}[-k]_q$ if $k<0$. Let $\t$ be a standard tableau and
suppose that~$i$ appears in row~$r$ and column~$c$ of~$\t$, where
$1\le i\le n$. The $q$--\textsf{content} of~$i$ in~$\t$ is
$c_\t(i)=[c-r]_q$.
Then, by \citeAM[Theorem~3.32]{M:Ulect},
$$m^\lambda_{\s\t}L_i=c_\t(i) m^\lambda_{\s\t}
               +\text{more dominant terms}.$$
Hence, $\{L_1,\dots,L_n\}$ is a family of JM--elements for $\H$.
Moreover, if $[1]_q[2]_q\dots[n]_q\ne0$ then a straightforward
induction shows that the JM--elements separate~$T(\Lambda)$;
see~\citeAM[Lemma~3.34]{M:Ulect}.
\end{Example}

\begin{Example}[Ariki--Koike algebras]\label{AK ex}
Fix integers $n,m\ge1$, an invertible element $q\in R$ and an
$m$--tuple $\mathbf u=(u_1\dots,u_m)\in R^m$. The \textsf{Ariki--Koike
algebra} $\H_{R,q,\mathbf u}$ is a deformation of the group algebra of
the complex reflection group of type $G(m,1,n)$; that is, the group
$(\Z/m\Z)\wr\Sym_n$. The Ariki--Koike algebras are generated by
elements $T_0,T_1,\dots,T_{n-1}$ subject to the relations
$(T_0-u_1)\dots(T_0-u_m)=0$, $(T_i-q)(T_i+1)=0$ for $1\le i<n$,
together with the braid relations of type~$B$. 

Let $\Lambda$ be the set of \textsf{$m$--multipartitions} of $n$; that
is, the set of $m$--tuples of partitions which sum to $n$. Then
$\Lambda$ is a poset ordered by dominance. If $\lambda\in\Lambda$ then
a \textsf{standard $\lambda$--tableau}  is an $m$--tuple of standard
tableau $\t=(\t^{(1)},\dots,\t^{(m)})$ which, collectively, contain the
numbers $1,\dots,n$ and where $\t^{(s)}$ has shape $\lambda^{(s)}$.
Let $T(\lambda)$ be the set of standard $\lambda$--tableaux ordered by
dominance~\citeAM[(3.11)]{DJM:cyc}. It is shown in \citeAM{DJM:cyc}
that the Ariki--Koike algebra has a cellular basis of the form
$\set{m^\lambda_{\s\t}|\lambda\in\Lambda\And\s,\t\in T(\lambda)}$.

For $i=1,\dots,n$ set $L_i=q^{1-i}T_{i-1}\dots T_1T_iT_1\dots T_{i-1}$.
These elements commute, are invariant under the $*$~involution
of $\H_{R,q,\mathbf u}$ and
$$m^\lambda_{\s\t}L_i=c_\t(i) m^\lambda_{\s\t}
            +\text{ more dominant terms},$$
where $c_\t(i)=u_sq^{c-r}$ if $i$ appears in row~$r$ and column~$c$ of
$\t^{(s)}$. All of these facts are proved in
\citeAM[\S3]{JM:cyc-Schaper}. Hence, $\{L_1,\dots,L_n\}$ is a family of
JM--elements for~$\H_{R,q,\mathbf u}$. In this case, if
$[1]_q\dots[n]_q\prod_{1\le i<j\le m}\prod_{|d|<n}(q^du_i-u_j)\ne0$
and $q\ne1$ then the JM--elements
separate~$T(\Lambda)$ by~\citeAM[Lemma~3.12]{JM:cyc-Schaper}.

There is an analogous family of JM--elements for the degenerate
Ariki--Koike algebras. See \citeAM[\S6]{AMR} for details.
\end{Example}

\begin{Example}[Schur algebras]\label{Schur ex}
Let $\Lambda$ be the set of partitions of $n$, ordered by dominance,
and for $\mu\in\Lambda$ let $\Sym_\mu$ be the corresponding Young
subgroup of $\Sym_n$ and set $m_\mu=\sum_{w\in\Sym_\mu}T_w\in\H$. Then
the \textsf{$q$--Schur algebra} is the endomorphism algebra
$$S_{R,q}(n)=\End_\H\( \bigoplus_{\mu\in\Lambda}m_\mu\H\).$$ 
For $\lambda\in\Lambda$ let $T(\lambda)$ be the set of semistandard
$\lambda$--tableaux, and let $T_\mu(\lambda)\subseteq T(\lambda)$ be
the set of semistandard $\lambda$--tableaux of type $\mu$;
see~\citeAM[\S4.1]{M:Ulect}. The main result of \citeAM{DJM:cyc} says
that $S_{R,q}(n)$ has a cellular basis
$\set{\varphi^\lambda_{\fS\fT}|\lambda\in\Lambda\And \fS,\fT\in
T(\lambda)}$ where the homomorphism $\varphi^\lambda_{\fS\fT}$ is
given by left multiplication by a sum of Murphy basis elements
$m^\lambda_{\s\t}\in\H$ which depend on~$\fS$ and~$\fT$.

Let $\mu=(\mu_1,\dots,\mu_k)$ be a partition in $\Lambda$.
For $i=1,\dots,k$ let $L^\mu_i$ be the endomorphism of  $m_\mu\H$
which is given by
$$L^\mu_i(m_\mu h) 
  =\sum_{j=\mu_1+\dots+\mu_{i-1}+1}^{\mu_1+\dots+\mu_i}L_jm_\mu h,$$
for all $h\in\H$. Here, $L_1,\dots,L_n$ are the JM--elements of the
Hecke algebra $\H$. We can consider $L^\mu_i$ to be an element of
$S_{R,q}(n)$. Using properties of the JM--elements of $\H$ it is easy to
check that the $L^\mu_i$ commute, that they are $*$--invariant and by
\citeAM[Theorem~3.16]{JM:Schaper} that
$$\varphi^\lambda_{\fS\fT}L^\mu_i=\begin{cases}
    c_\fT(i)\varphi^\lambda_{\fS\fT}
          +\text{more dominant terms},&\If\fT\in T_\mu(\lambda),\\
    0,&\Otherwise.
\end{cases}$$
Here $c_\fT(i)$ is the sum of the $q$--contents of the nodes in $\fT$
labelled by~$i$~\citeAM[\S5.1]{M:Ulect}. Hence
$\set{L^\mu_i|\mu\in\Lambda}$ is a family of JM--elements for
$S_{R,q}(n)$. If $[1]_q\dots[n]_q\ne0$ then the JM--elements separate
$T(\Lambda)$; see \citeAM[Lemma~5.4]{M:Ulect}.

More generally, the $q$--Schur algebras $S_{R,q}(n,r)$ of type $A$ and the
cyclotomic $q$--Schur algebras both have a family of JM--elements; see
\citeAM{JM:Schaper,JM:cyc-Schaper} for details.
\end{Example}

\begin{Example}[Birman--Murakami--Wenzl algebras]\label{BMW ex}
Let $r$ and $q$ be invertible indeterminates over $R$ and let $n\ge1$ an
integer.  Let $\B$ be the Birman--Murakami--Wenzl algebra, or
\textsf{BMW algebra}. The BMW algebra is generated by elements
$T_1,\dots,T_{n-1}$ which satisfy the relations
$(T_i-q)(T_i+q^{-1})(T_i-r^{-1})=0$, the braid relations of type~$A$,
and the relations $E_iT^{\pm1}_{i\pm1}E_i=r^{\pm1} E_i$ and 
$E_iT_i=T_iE_i=r^{-1}E_i$, where $E_i=1-\frac{T_i-T_i^{-1}}{q-q^{-1}}$;
see \citeAM{Enyang:ss,LeducRam}. 

The BMW algebra $\B$ is a deformation of the Brauer algebra.
Indeed, both the Brauer and BMW algebras have a natural diagram basis
indexed by the set of $n$--Brauer diagrams; that is, graphs with
vertex set $\{1,\dots,n,\bar 1,\dots,\bar n\}$ such that each vertex
lies on a unique edge. For more details see \citeAM{HalRam:basic}.

Let $\lambda$ be a partition of $n-2k$, where $0\le k\le\floor n2$.
An $n$--updown $\lambda$--tableau $\t$ is an $n$--tuple
$\t=(\t_1,\dots,\t_n)$ of partitions such that $\t_1=(1)$,
$\t_n=\lambda$ and $|\t_{i}|=|\t_{i-1}|\pm1$, for $2\le i\le n$. (Here
$|\t_i|$ is the sum of the parts of the partition $\t_i$.)

Let $\Lambda$ be the set of partitions of $n-2k$, for $0\le k\le\floor
n2$ ordered again by dominance. For $\lambda\in\Lambda$ let
$T(\lambda)$ be the set of $n$--updown tableaux.
Enyang~\citeAM[Theorem~4.8 and \S5]{Enyang:ss} has given an algorithm
for constructing a cellular basis of $\B$ of the form
$\set{m^\lambda_{\s\t}|\lambda\in\Lambda\And\s,\t\in T(\lambda)}$.
Enyang actually constructs a basis for each cell module of $\B$
which is ``compatible'' with restriction, however, his arguments give
a  new cellular basis $\{m_{\s\t}^\lambda\}$ for $\B$ which is indexed
by pairs of $n$--updown $\lambda$--tableaux for $\lambda\in\Lambda$.

Following \citeAM[Cor.~1.6]{LeducRam} set $L_1=1$ and define
$L_{i+1}=T_iL_iT_i$, for $i=2,\dots,n$. These elements are invariant
under the $*$~involution of $\B$ and
Enyang~\citeAM[\S6]{Enyang:ss} has shown that $L_1.\dots,L_n$ 
commute and that
$$m^\lambda_{\s\t}L_i=c_\t(i)m^\lambda_{\s\t}
                    +\text{ more dominant terms},$$
where $c_\t(i)=q^{2(c-r)}$ if $[\t_i]=[\t_{i-1}]\cup\{(r,c)\}$ and
$c_\t(i)=r^{-2}q^{2(r-c)}$ if $[\t_i]=[\t_{i-1}]\setminus\{(r,c)\}$.
Hence, $L_1,\dots,L_n$ is a family of JM--elements for $\B$. When
$R=\Z[r^{\pm1},q^{\pm1}]$ the JM--elements separate $T(\Lambda)$.

The BMW algebras include the Brauer algebras essentially as a special
case. Indeed, it follows from Enyang's work~\citeAM[\S8--9]{Enyang:ss}
that the Brauer algebras have a family of JM--elements which
separate~$T(\Lambda)$. 

Rui and Si~\citeAM{RuiSi:BrauerDet} have recently computed the Gram
determinants of the irreducible modules of the Brauer algebras
in the semisimple case.
\end{Example}

It should be possible to find JM--elements for other cellular algebras
such as the partition algebras and the cyclotomic Nazarov--Wenzl
algebras~\citeAM{AMR}.

\section{The separated case}

Throughout this section we assume that $A$ is a cellular algebra with
a family of JM--elements which separate $T(\Lambda)$ over~$R$. By
Corollary~\ref{semisimple} this implies that $A_K$ is a split
semisimple algebra. 

For $i=1,\dots,M$ let $\CC(i)=\set{c_\t(i)|\t\in T(\Lambda)}$. Thus,
$\CC(i)$ is the set of possible contents that the elements of
$T(\Lambda)$ can take at~$i$.

We can now make the key definition of this paper.

\begin{Defn}\label{F_t}Suppose that $\s,\t\in T(\lambda)$, for some
$\lambda\in\Lambda$ and define 
$$F_\t=\Prod_{i=1}^M\prod_{\substack{c\in\CC(i)\\c\ne c_\t(i)}}
                    \frac{L_i-c}{c_\t(i)-c}.$$
Thus, $F_\t\in A_K$. Define 
$f^\lambda_{\s\t}=F_\s a^\lambda_{\s\t} F_\t\in A_K$.
\end{Defn}

\begin{Remark} 
Rather than working over $K$ we could instead work over a ring $R'$ in
which the elements
$\set{c_\s(i)-c_\t(i)|\s\ne\t\in T(\Lambda)\And 1\le i\le M}$ are
invertible. All of the results below, except those concerned with the
irreducibe $A_K$--modules or with the semisimplicity of $A_K$, are
valid over~$R'$. However, there seems to be no real advantage to
working over~$R'$ in this section. In section~4 we work over a similar
ring when studying the non-separated case.
\end{Remark}

We extend the dominance order $\gdom$ on $T(\Lambda)$ to
$\coprod_{\lambda\in\Lambda}T(\lambda)\times T(\lambda)$ by declaring
that $(\s,\t)\gdom(\u,\v)$ if $\s\gedom\u$, $\t\gedom\v$ and
$(\s,\t)\ne(\u,\v)$.

We now begin to apply our definitions. The first step is easy.

\begin{Lemma}\label{f basis}
Assume that $A$ has a family of JM--elements which separate~$T(\Lambda)$. 
\begin{enumerate}
\item Suppose that $\s,\t\in T(\lambda)$.  Then there exist scalars
$b_{\u\v}\in K$ such that
$$f^\lambda_{\s\t}=a^\lambda_{\s\t}
+\sum_{\substack{\u,\v\in T(\mu),\mu\in\Lambda\\(\u,\v)\gdom(\s,\t)}}
              b_{\u\v}a^\mu_{\u\v}.$$
\item $\set{f^\lambda_{\s\t}|\s,\t\in T(\lambda)\text{ for some
}\lambda\in\Lambda}$ is a basis of $A_K$.
\item Suppose that $\s,\t\in T(\lambda)$. Then and
$(f^\lambda_{\s\t})^*=f^\lambda_{\t\s}$.
\end{enumerate}
\end{Lemma}

\begin{proof}
By the definition of the JM--elements (\ref{JM}), for any $i$ and any
$c\in\CC(i)$ with $c\ne c_\t(i)$ we have
$$a^\lambda_{\s\t}\frac{L_i-c}{c_\t(i)-c}
         \equiv a^\lambda_{\s\t}+\sum_{\v\gdom\t}b_\v a^\lambda_{\s\v}
           \pmod{A_K^\lambda}.$$
By (\ref{JM*}) this is still true if we act on
$a^\lambda_{\s\t}$ with $L_i$ from the left. These two facts imply
part~(a).  Note that part (a) says that the transition matrix between
the two bases $\{a^\lambda_{\s\t}\}$ and $\{f^\lambda_{\s\t}\}$ of
$A_K$ is unitriangular (when the rows and columns are suitably
ordered). Hence,~(b) follows.  Part (c) follows because, by definition,
$(a^\lambda_{\s\t})^*=a^\lambda_{\t\s}$ and $L_i^*=L_i$, so that
$F_\t^*=F_\t$ and
$(f^\lambda_{\s\t})^*=F_\t a^\lambda_{\t\s}F_\s=f^\lambda_{\t\s}$.
\end{proof}

Given $\s,\t\in T(\Lambda)$ let $\delta_{\s\t}$ be the Kronecker
delta; that is, $\delta_{\s\t}=1$ if $\s=\t$ and $\delta_{\s\t}=0$,
otherwise.

\begin{Prop}\label{vanishing}
Suppose that $\s,\t\in T(\lambda)$, for some
$\lambda\in\Lambda$, that $\u\in T(\Lambda)$ and fix~$i$ with 
$1\le i\le M$. Then
\begin{multicols}{2}
\begin{enumerate}
\item $f^\lambda_{\s\t}L_i=c_\t(i)f^\lambda_{\s\t}$,
\item $f^\lambda_{\s\t}F_\u=\delta_{\t\u}f^\lambda_{\s\u}$,
\item $L_if^\lambda_{\s\t}=c_\s(i)f^\lambda_{\s\t}$,
\item $F_\u f^\lambda_{\s\t}=\delta_{\u\s}f^\lambda_{\u\t}$.
\end{enumerate}
\end{multicols}
\end{Prop}

\begin{proof}
Notice that statements (a) and (c) are equivalent by applying the
$*$~involution. Similarly, (b) and~(d) are equivalent.  Thus, it is
enough to show that~(a) and~(b) hold.  Rather than proving this
directly we take a slight detour.

Let $N=|T(\Lambda)|$ and fix $\v=\v_1\in T(\mu)$ with $\v\gdom\t$.  We
claim that $a^{\mu}_{\u\v}F_\t^N=0$, for all $\u\in T(\mu)$. By the
separation condition (\ref{separation}), there exists an integer $j_1$
with $c_\t(j_1)\ne c_\v(j_1)$. Therefore, by (\ref{JM}),
$a^{\mu}_{\u\v}(L_{j_1}-c_\v(j_1))$ is a linear combination of terms
$a^\nu_{\w\x}$, where $\x\gdom\v\gdom\t$.  However,
$(L_{j_1}-c_\v(j_1))$ is a factor of $F_\t$, so~$a^{\mu}_{\u\v}F_\t$
is a linear combination of terms of the form $a^\nu_{\w\x}$ where
$\x\gdom\v\gdom\t$. Let $\v_2\in T(\mu_2)$ be minimal such that
$a^{\mu_2}_{\u_2\v_2}$ appears with non--zero coefficient in
$a^{\mu}_{\u\v}F_\t$, for some $\u_2\in T(\mu_2)$. Then
$\v_2\gdom\v_1\gdom\t$, so there exists an integer $j_2$ such that
$c_\t(j_2)\ne c_{\v_2}(j_2)$. Consequently, $(L_{j_2}-c_{\v_2}(j_2))$
is a factor of $F_\t$,  so $a^\mu_{\u\v}F_\t^2$ is a linear
combination of terms of the form $a^\nu_{\w\x}$, where
$\x\gdom\v_2\gdom\v_1\gdom\t$.  Continuing in this way proves the claim.

For any $\s,\t\in T(\lambda)$ let 
$f_{\s\t}'=F_\s^N a^\lambda_{\s\t}F_\t^N$. Fix $j$ with $1\le j\le M$. Then, 
because the JM--elements commute,
$$ f_{\s\t}'L_j=F_\s^N a^\lambda_{\s\t}F_\t^NL_j
      =F_\s^N a^\lambda_{\s\t}L_jF_\t^N    
      =F_\s^N\(c_\t(i)a^\lambda_{\s\t}+x\)F_\t^N,$$
where $x$ is a linear combination of terms of the form $a^\mu_{\u\v}$
with $\v\gdom\t$ and $\u,\v\in T(\mu)$ for some $\mu\in\Lambda$.
However, by the last paragraph $xF_\t^N=0$, so this implies that
$f_{\s\t}'L_j=c_\t(j)f_{\s\t}'$.  Consequently, every factor of $F_\t$
fixes $f_{\s\t}'$, so $f_{\s\t}'=f_{\s\t}'F_\t$.  Moreover, if
$\u\ne\t$ then we can find $j$ such that $c_\t(j)\ne c_\u(j)$ by the
separation condition, so that $f_{\s\t}'F_\u=0$ since $(L_j-c_\u(j))$
is a factor of $F_\u$. As $F_\u f_{\s\t}'=(f_{\t\s}'F_\u)^*$, 
we have shown that 
\begin{equation}\label{vanishing'}
    F_\u f_{\s\t}'F_\v=\delta_{\u\s}\delta_{\t\v}f_{\s\t}',
\end{equation}
for any $\u,\v\in T(\Lambda)$.

We are now almost done. By the argument of Lemma~\ref{f basis}(a) we 
know that 
$$f_{\s\t}'=a^\lambda_{\s\t}
    +\sum_{\substack{\u,\v\in T(\mu)\\(\u,\v)\gdom(\s,\t)}}
             s_{\u\v}a^\mu_{\u\v},$$
for some $s_{\u\v}\in K$. Inverting this equation we can write
$$a^\lambda_{\s\t}=f_{\s\t}'
    +\sum_{\substack{\u,\v\in T(\mu)\\(\u,\v)\gdom(\s,\t)}}
             s_{\u\v}'f_{\u\v}',$$
for some $s_{\u\v}'\in K$. Therefore,
$$f^\lambda_{\s\t}
    =F_\s a^\lambda_{\s\t}F_\t
    =F_\s\(f_{\s\t}'
    +\sum_{\substack{\u,\v\in T(\mu)\\(\u,\v)\gdom(\s,\t)}}
             s_{\u\v}'f_{\u\v}'\) F_\t
    =F_\s f_{\s\t}' F_\t=f_{\s\t}',
$$
where the last two equalities follow from (\ref{vanishing'}). That
is, $f^\lambda_{\s\t}=f_{\s\t}'$. We now have that
$$f^\lambda_{\s\t}L_i=f_{\s\t}'L_i=c_\t(i)f_{\s\t}'
          =c_\t(i)f^\lambda_{\s\t},$$
proving (a). Finally, if $\u\in T(\Lambda)$ then
$$f^\lambda_{\s\t}F_\u=f_{\s\t}'F_\u=\delta_{\t\u}f_{\s\t}'
          =\delta_{\t\u}f^\lambda_{\s\t},$$
proving (b). (In fact, (b) also follows from (a) and the separation
condition.)
\end{proof}

\begin{Remark}
    The proof of the Proposition~\ref{vanishing} is the only place where we
    explicitly invoke the separation condition. All of the
    results which follow rely on this key result. It is also worth
    noting the proof of Proposition~\ref{vanishing} relies on
    the two assumptions that the $L_1,\dots,L_M$ commute and that $L_i^*=L_i$,
    for $\le i\le M$. The commutivity of the JM--elements is essential
    in proving (\ref{vanishing}). If we did not assume that
    $L_i\ne L_i^*$ then we could define $f_{\s\t}=F_\s^* a_{\s\t}^\lambda F_\t$. If we did this then 
    in order to prove that $f_{s\t}'F_\u=0$ we would assume that the~$L_i^*$ act from the right on the basis $\{a_{\v}^\mu\}$
    in essentialy the same way as the $L_j$ do. We note that neither
    of these assumptions appear in Soriano's treatment in the appendix.
\end{Remark}

\begin{Theorem}\label{orthogonal}
Suppose that the JM--elements separate $T(\Lambda)$ over $R$. Let 
$\s,\t\in T(\lambda)$ and $\u,\v\in T(\mu)$,
for some $\lambda,\mu\in\Lambda$. Then there exist scalars
$\set{\gamma_\t\in K|\t\in T(\Lambda)}$ such that
$$f^\lambda_{\s\t}f^\mu_{\u\v}=\begin{cases}
   \gamma_\t f^\lambda_{\s\v},&\text{if }\lambda=\mu\And\t=\u,\\
           0,&\text{otherwise.}
\end{cases}$$
In particular, $\gamma_\t$ depends only on $\t\in T(\Lambda)$ and
$\set{f^\lambda_{\s\t}|\s,\t\in T(\lambda)\And\lambda\in\Lambda}$
is a cellular basis of $A_K$.
\end{Theorem}

\begin{proof}Using the definitions,
$f^\lambda_{\s\t}f^\mu_{\u\v}=f^\lambda_{\s\t} F_\u a^\mu_{\u\v}F_\v$.
So $f^\lambda_{\s\t}f^\mu_{\u\v}\ne0$ only if $\u=\t$ by
Proposition~\ref{vanishing}(b). 

Now suppose that $\u=\t$ (so that $\mu=\lambda$). Using 
Lemma~\ref{f basis}, we can write
$f^\lambda_{\s\t}f^\lambda_{\t\v}=\sum_{\w,\x}r_{\w\x}f^\mu_{\w\x}$,
where $r_{\w\x}\in R$ and the sum is over pairs $\w,\x\in T(\mu)$, for
some $\mu\in\Lambda$. Hence, by parts (b) and (d) of
Proposition~\ref{vanishing}
\begin{equation}\label{almost}
f^\lambda_{\s\t}f^\lambda_{\t\v}
       =F_\s f^\lambda_{\s\t}f^\lambda_{\t\v} F_\v
       =\sum_{\substack{\mu\in\Lambda\\\w,\x\in T(\mu)}}
               r_{\w\x}F_\s f^\mu_{\w\x}F_\v
       =r_{\s\v}f^\lambda_{\s\v}.
\end{equation}
Thus, it remains to show that scalar $r_{\s\v}$ is independent of
$\s,\v\in T(\lambda)$. Using Lemma~\ref{f basis} to compute directly,
there exist scalars $b_{\w\x},c_{\y\z},r_{\w\x} \in K$ such that
\begin{align*}
f^\lambda_{\s\t}f^\lambda_{\t\v}
    &\equiv
      \(a^\lambda_{\s\t}
      +\sum_{\substack{\w,\x\in T(\lambda)\\(\w,\x)\gdom(\s,\t)}}
                b_{\w\x} a^\lambda_{\w\x}\)
      \(a^\lambda_{\t\v}
      +\sum_{\substack{\y,\z\in T(\lambda)\\(\y,\z)\gdom(\t,\v)}}
               c_{\y\z}a^\lambda_{\y\z}\) \mod{A_K^\lambda}\\
   &\equiv \(\<a^\lambda_\t,a^\lambda_\t\>_\lambda
   +\sum_{\substack{\u\in T(\lambda)\\\u\gdom\t}}
        (b_{\s\u}+c_{\u\v})\<a^\lambda_\u,a^\lambda_\t\>_\lambda
   +\sum_{\substack{\x,\y\in T(\lambda)\\\x,\y\gdom\t}}
	b_{\s\x}c_{\y\v}\<a^\lambda_\x,a^\lambda_\y\>_\lambda\)
	    a^\lambda_{\s\v}\\
   &\hspace*{20mm} +\sum_{\substack{\w,\x\in T(\lambda)\\(\w,\x)\gdom(\s,\v)}}
                  r_{\w\x}a^\lambda_{\w\x}\quad \pmod{A_K^\lambda}.
\end{align*}
The inner products in the last  equation come from applying (\ref{inner product}).
(For typographical convenience we also use the fact that the form is
symmetric in the sum over~$\u$.) That is, there exists a scalar
$\gamma\in A$, which does not depend on $\s$ or on $\v$, such that
$f^\lambda_{\s\t}f^\lambda_{\t\v}=\gamma a^\lambda_{\s\v}$ plus a linear
combination of more dominant terms. By Lemma~\ref{f basis}(b) and
(\ref{almost}), the coefficient of $f^\lambda_{\s\v}$ 
in~$f^\lambda_{\s\t}f^\lambda_{\t\v}$ is equal to the coefficient of
$a^\lambda_{\s\v}$ in~$f^\lambda_{\s\t}f^\lambda_{\t\v}$, so this
completes the proof.
\end{proof}

We call 
$\set{f^\lambda_{\s\t}|\s,\t\in T(\lambda)\And\lambda\in\Lambda}$ 
the \textsf{seminormal basis} of~$A$. This terminology is justified by
Remark~\ref{orthogonality} below.

\begin{Cor}\label{non-zero gammas}
Suppose that $A_K$ is a cellular algebra with a family of
JM--elements which separate $T(\Lambda)$. Then $\gamma_\t\ne0$, for all
$\t\in T(\Lambda)$.
\end{Cor}

\begin{proof}
Suppose by way of contradiction that $\gamma_\t=0$, for some $\t\in
T(\lambda)$ and $\lambda\in\Lambda$. Then, by
Theorem~\ref{orthogonal},
$f^\lambda_{\t\t}f^\mu_{\u\v}=0=f^\mu_{\u\v}f^\lambda_{\t\t}$, for all
$\u,\v\in T(\mu)$, $\mu\in\Lambda$.  Therefore, $Kf^\lambda_{\t\t}$
is a one dimensional nilpotent ideal of~$A_K$, so~$A_K$ is not
semisimple. This contradicts Corollary~\ref{semisimple}, so
we must have $\gamma_\t\ne0$ for all $\t\in T(\Lambda)$.
\end{proof}

Next, we use the basis $\{f^\lambda_{\s\t}\}$ to identify the cell
modules of $A$ as submodules of~$A$.

\begin{Cor}\label{cyclic}
Suppose that $\lambda\in\Lambda$ and fix $\s,\t\in T(\lambda)$. Then
$$C(\lambda)\cong f^\lambda_{\s\t}A_K
             = \Span\set{f^\lambda_{\s\v}|\v\in T(\lambda)}.$$
\end{Cor}

\begin{proof}
As $f^\mu_{\u\v}=F_\u a^\mu_{\u\v}F_\v$, for $\u,\v\in T(\mu)$, the cell
modules for the cellular bases $\{a^\lambda_{\u\v}\}$ and $\{f^\lambda_{\u\v}\}$ 
of $A_K$ coincide. Therefore, $C(\lambda)$ is isomorphic to the
$A_K$--module $C(\lambda)'$ which is spanned by the elements
$\set{f^\lambda_{\s\u}+A_K^\lambda|\u\in T(\lambda)}$.

On other hand, if $\u,\v\in T(\mu)$, for $\mu\in\Lambda$, then
$f^\lambda_{\s\t}f^\mu_{\u\v}=\delta_{\t\u}\gamma_\t f^\lambda_{\s\v}$
by Theorem~\ref{orthogonal}.  Now $\gamma_\t\ne0$, by
Corollary~\ref{non-zero gammas}, so $\set{f^\lambda_{\s\v}|\v\in
T(\lambda)}$ is a basis of $f^\lambda_{\s\t}A_K$.

Finally, by Theorem~\ref{orthogonal}  we have that 
$f^\lambda_{\s\t}A_K\cong C(\lambda)'$, where the
isomorphism is the linear extension of the map
$f^\lambda_{\s\v}\mapsto f^\lambda_{\s\v}+A_K^\lambda$,
for $\v\in T(\lambda)$. Hence, 
$C(\lambda)\cong C(\lambda)'\cong f^\lambda_{\s\t}A_K$, as required.
\end{proof}

Recall that $\rad C(\lambda)$ is the radical of the bilinear form on
$C(\lambda)$ and that $D(\lambda)=C(\lambda)/\rad C(\lambda)$.

Using Corollary~\ref{cyclic} and Theorem~\ref{orthogonal}, the basis
$\{f^\lambda_{\s\t}\}$ gives an explicit decomposition of $A_K$ into a
direct sum of cell modules. Abstractly this also follows from
Corollary~\ref{semisimple} and the general theory of cellular algebras
because a cellular algebra is semisimple if and only if
$C(\lambda)=D(\lambda)$, for all $\lambda\in\Lambda$;
see~\citeAM[Theorem~3.4]{GL}.

\begin{Cor}\label{cell decomp}
Suppose that $A_K$ is a cellular algebra with a family of
JM--elements which separate $T(\Lambda)$. Then $C(\lambda)=D(\lambda)$, for
all $\lambda\in\Lambda$, and
$$A_K\cong\bigoplus_{\lambda\in\Lambda} 
                C(\lambda)^{\oplus|T(\lambda)|}.$$
\end{Cor}

Fix $\s\in T(\lambda)$ and, for notational convenience, set
$f^\lambda_\t=f^\lambda_{\s\t}$ so that $C(\lambda)$ has basis
$\set{f^\lambda_\t|\t\in T(\lambda)}$ by Corollary~\ref{cyclic}.
Note that $f^\lambda_\t=a^\lambda_\t+\sum_{\v\gdom\t}b_\v
a^\lambda_\v$, for some $b_\v\in K$, by Lemma~\ref{f basis}(a).

For $\lambda\in\Lambda$ let
$G(\lambda)=\det\big(\<a^\lambda_\s,a^\lambda_\t\>_\lambda%
	   \big)_{\s,\t\in T(\lambda)}$ 
be the Gram determinant of the bilinear form $\<\ ,\ \>_\lambda$ on
the cell module $C(\lambda)$. Note that $G(\lambda)$ is well--defined
only up to multiplication by $\pm1$ as we have not specified an
ordering on the rows and columns of the Gram matrix.

\begin{Theorem}\label{Gram det}
Suppose that $A_K$ is a cellular algebra with a family of
JM--elements which separate $T(\Lambda)$. Let $\lambda\in\Lambda$ and suppose
that $\s,\t\in T(\lambda)$. Then
$$\<f^\lambda_\s,f^\lambda_\t\>_\lambda
          =\<a^\lambda_\s,f^\lambda_\t\>_\lambda
          =\begin{cases}\gamma_\t,&\If\s=\t,\\
              0,&\Otherwise.
\end{cases}$$
Consequently, $G(\lambda)=\Prod_{\t\in T(\lambda)}\gamma_\t$.
\end{Theorem}

\begin{proof}By Theorem~\ref{orthogonal}, $\{f^\lambda_{\s\t}\}$ is a
cellular basis of  $A_K$ and, by Corollary~\ref{cyclic}, we may take
$\set{f^\lambda_\t|\t\in T(\lambda)}$ to be a basis of $C(\lambda)$.
By Theorem~\ref{orthogonal} again,
$f^\lambda_{\u\s}f^\lambda_{\t\v}=\delta_{\s\t}\gamma_\t f^\lambda_{\u\v}$,
so that $\<f^\lambda_\s,f^\lambda_\t\>_\lambda=\delta_{\s\t}\gamma_\t$
by Corollary~\ref{cyclic} and the definition of the inner product on
$C(\lambda)$. Using Proposition~\ref{vanishing}(b) and the
associativity of the inner product on $C(\lambda)$, we see that
$$\<a^\lambda_\s,f^\lambda_\t\>_\lambda
     =\<a^\lambda_\s,f^\lambda_\t F_\t\>_\lambda
     =\<a^\lambda_\s F_\t^*,f^\lambda_\t\>_\lambda
     =\<a^\lambda_\s F_\t,f^\lambda_\t\>_\lambda
     =\<f^\lambda_\s,f^\lambda_\t\>_\lambda.$$
So we have proved the first claim in the statement of the Theorem.

Finally, the transition matrix between the two bases
$\{a^\lambda_\t\}$ and $\{f^\lambda_\t\}$ of~$C(\lambda)$ is
unitriangular (when suitably ordered), so we have that
$$G(\lambda)=\det\big(\<a^\lambda_\s,a^\lambda_\t\>_\lambda\big)
           =\det\big(\<f^\lambda_\s,f^\lambda_\t\>_\lambda\big)
           =\prod_{\t\in T(\lambda)}\gamma_\t,$$
as required.
\end{proof}

\begin{Remark}\label{orthogonality}
Extending the bilinear forms $\<\ , \>_\lambda$ to the whole of $A_K$
(using Corollary~\ref{cell decomp}), we see that the seminormal basis
$\{f^\lambda_{\s\t}\}$ is an orthogonal basis of $A_K$ with respect to
this form.
\end{Remark}

In principle, we can use Theorem~\ref{Gram det} to compute the Gram
determinants of the cell modules of any cellular algebra $A$ which has
a separable family of JM--elements. In practice, of course, we need to
find formulae for the structure constants $\set{\gamma_\t|\t\in
T(\lambda)}$ of the basis $\{f^\lambda_{\s\t}\}$. In all known
examples, explicit formulae for~$\gamma_\t$ can be determined
inductively once the actions of the generators of $A$ on the
seminormal basis have been determined. In turn, the action of $A$ on
its seminormal basis is determined by its action on the original
cellular basis $\{a^\lambda_{\s\t}\}$. In effect, Theorem~\ref{Gram
det} gives an effective recipe for computing the Gram determinants of
the cell modules of~$A$.

By definition the scalars $\gamma_\t$ are elements of the field~$K$,
for $\t\in T(\lambda)$. Surprisingly, their product must belong to
$R$.

\begin{Cor} Suppose that $\lambda\in\Lambda$. Then $\Prod_{\t\in
T(\lambda)}\gamma_\t\in R$.
\end{Cor}

\begin{proof}By definition, the inner products 
$\<a^\lambda_\s,a^\lambda_\t\>_\lambda$ all belong to $R$, so 
$G(\lambda)\in R$. The result now follows from Theorem~\ref{Gram det}.
\end{proof}

As $G(\lambda)\ne0$ by Theorem~\ref{Gram det} and
Corollary~\ref{non-zero gammas}, it follows that each cell module is
irreducible.

\begin{Cor}\label{irred cells}
Suppose that $\lambda\in\Lambda$. Then the cell module
    $C(\lambda)=D(\lambda)$ is irreducible.
\end{Cor}

We close this section by describing the primitive idempotents in~$A_K$.

\begin{Theorem}\label{idempotents}
Suppose that $A_K$ is a cellular algebra with a family of
JM--elements which separate $T(\Lambda)$. Then
\begin{enumerate}
\item If $\t\in T(\lambda)$ and $\lambda\in\Lambda$ then 
$F_\t=\frac1{\gamma_\t}f^\lambda_{\t\t}$ and $F_\t$ is a primitive 
idempotent in~$A_K$.
\item If $\lambda\in\Lambda$ then $F_\lambda=\sum_{\t\in
T(\lambda)}F_\t$ is a primitive central idempotent in $A_K$.
\item $\set{F_\t|\t\in T(\Lambda)}$ and $\set{F_\lambda|\lambda\in\Lambda}$ 
    are complete sets of pairwise orthogonal idempotents in $A_K$; 
    in particular,
$$1_{A_K}=\sum_{\lambda\in\Lambda}F_\lambda
     =\sum_{\t\in T(\Lambda)}F_\t.$$
\end{enumerate}
\end{Theorem}

\begin{proof}By Corollary~\ref{non-zero gammas}, $\gamma_\t\ne0$ for
all $\t\in T(\lambda)$, so the statement of the Theorem makes
sense. Furthermore, $\frac1{\gamma_\t}f^\lambda_{\t\t}$
is an idempotent by Theorem~\ref{orthogonal}. By 
Corollary~\ref{irred cells} the cell module $C(\lambda)$ is
irreducible and by Corollary~\ref{cyclic},
$C(\lambda)\cong f^\lambda_{\t\t}A_K=F_\t A_K$. Hence, $F_\t$ is a
primitive idempotent. 

To complete the proof of (a) we still need to show that
$F_\t=\frac1{\gamma_\t}f^\lambda_{\t\t}$. By Theorem~\ref{orthogonal}
we can write
$F_\t=\sum_{\nu\in\Lambda}\sum_{\x,\y\in T(\nu)} r_{\x\y}f^\nu_{\x\y}$, 
for some $r_{\x\y}\in K$. Suppose that $\u,\v\in T(\mu)$, for some
$\mu\in\Lambda$. Then, by Proposition~\ref{vanishing} and
Theorem~\ref{orthogonal},
$$\delta_{\v\t}f^\mu_{\u\v}
    =f^\mu_{\u\v}F_\t 
    =\sum_{\nu\in\Lambda}\sum_{\x,\y\in T(\mu)}
      r_{\x\y}f^\mu_{\u\v}f^\nu_{\x\y}
      =\sum_{\y\in T(\mu)} r_{\v\y}\gamma_\v f^\mu_{\u\y}.$$
By Corollary~\ref{non-zero gammas} $\gamma_\v\ne0$, so comparing
both sides of this equation shows that
$$r_{\v\y}=\begin{cases}\frac1{\gamma_\t},&\If \v=\t=\y,\\
          0,&\Otherwise.
\end{cases}$$
As $\v$ is arbitrary we have 
$F_\t=\frac1{\gamma_\t}f^\lambda_{\t\t}$, as claimed.

This completes the proof of~(a). Parts (b) and (c) now follow from (a)
and the multiplication formula in Theorem~\ref{orthogonal}.
\end{proof}

\begin{Cor}\label{minimum poly}
Suppose that $A_K$ is a cellular algebra with a family of
JM--elements which separate $T(\Lambda)$. Then
$$L_i=\sum_{\t\in T(\Lambda)}c_\t(i)F_\t$$
and $\prod_{c\in\CC(i)}(L_i-c)$ is the minimum polynomial for $L_i$
acting on $A_K$.
\end{Cor}

\begin{proof}By part (c) of Theorem~\ref{idempotents},
$$L_i = L_i\sum_{\t\in T(\Lambda)} F_\t
      = \sum_{\t\in T(\Lambda)} L_iF_\t
      = \sum_{\t\in T(\Lambda)} c_\t(i)F_\t,$$
where the last equality follows from Proposition~\ref{vanishing}(c).

For the second claim, observe that 
$\prod_{c\in\CC(i)} (L_i-c)\cdot f^\lambda_{\s\t}=0$ by
Proposition~\ref{vanishing}(c), for all $\lambda\in\Lambda$ and all
$\s,\t\in T(\lambda)$. If we omit the factor $(L_i-d)$, for some
$d\in\CC(i)$, then we can find an $\s\in T(\mu)$, for some $\mu$, such
that $c_\s(i)=d$ so that $\prod_{c\ne d}(L_i-c)F_\s\ne0$. Hence,
$\prod_{c\in\CC(i)}(L_i-c)$ is the minimum polynomial for the action
of $L_i$ on $A_K$.
\end{proof}

The examples at the end of section~2 show that the number of
JM--elements is not uniquely determined. Nonetheless, we are able to
characterize the subalgebra of~$A_K$ which they generate.

\begin{Cor}Suppose that $A_K$ is a cellular algebra with a family of
JM--elements which separate $T(\Lambda)$. Then
$\{L_1,\dots,L_M\}$ generate a maximal abelian subalgebra of~$A_K$.
\end{Cor}

\begin{proof}As the JM--elements commute, by definition, the
subalgebra $\L_K$ of $A_K$ which they generate is certainly abelian. By
Theorem~\ref{idempotents} and Corollary~\ref{minimum poly}, $\L_K$ is the
subalgebra of~$A$ spanned by the primitive idempotents 
$\set{F_\t|\t\in T(\Lambda)}$. As the primitive idempotents of $A_K$
span a maximal abelian subalgebra of $A_K$, we are done.
\end{proof}

\section{The non--separated case}
Up until now we have considered those cellular algebras $A_K$ which
have a family of JM--elements which separate $T(\Lambda)$. By
Corollary~\ref{semisimple} the separation condition forces $A_K$ to be
semisimple. In this section we still assume that $A=A_R$ has a family
of JM--elements which separate $T(\Lambda)$ over~$R$ but rather than
studying the semisimple algebra $A_K$ we extend the previous
constructions to non--separated algebras over a field.

In this section let $R$ be a discrete valuation ring with maximal
ideal~$\pi$. We assume that~$A_R$ has a family of JM--elements which
separate $T(\Lambda)$ over~$R$. 

Let $K$ be the field of fractions of $R$. Then $A_K$ is semisimple by
Corollary~\ref{semisimple} and all of the results of the previous
section apply to~$A_K$. Let $k=R/\pi$ be the residue field of~$K$.
Then $A_k=A\otimes_R k$ is a cellular algebra with cellular basis
given by the image of the cellular basis of~$A$ in~$A_k$. We abuse
notation and write $\{a^\lambda_{\s\t}\}$ for the cellular bases of
all three algebras $A=A_R$, $A_K$ and~$A_k$. It should always be clear
from the context which algebra these elements belong to at any given
time.

In general, the JM--elements will not separate $T(\Lambda)$ over~$k$,
so the arguments of the previous section do not necessarily apply to
the algebra~$A_k$. 

If $r\in R$ let $\bar r=r+\pi$ be its image in $k=R/\pi$. More
generally, if $a=\sum r_{\s\t}a^\lambda_{\s\t}\in A_R$ then we set
$\bar a=\sum \bar{r_{\s\t}}a^\lambda_{\s\t}\in A_k$. 

The final assumption that we make is that $c-c'$ is invertible in $R$
whenever $\bar c\ne\bar c'$ and~$c,c'\in\CC=\bigcup_{i=1}^M\CC(i)$.

If $1\le i\le M$ and $\t\in T(\lambda)$ define the \textsf{residue} of
$i$ at $\t$ to be $r_\t(i)=\bar{c_\t(i)}$. 
By (\ref{JM}) the action of the JM--elements on $A_k$ is given by
$$ a^\lambda_{\s\t}L_i
   \equiv r_\t(i)a^\lambda_{\s\t}
        +\sum_{\v\gdom\t}r_{\t\v}a^\lambda_{\s\v} \pmod{A_k^\lambda},
$$
where $r_{\t\v}\in k$ (and otherwise the notation is as in (\ref{JM})).
There is an analogous formula for the action of $L_i$ on
$a^\lambda_{\s\t}$ from the left.

We use residues modulo $\pi$ to define equivalence relations on
$T(\Lambda)$ and on $\Lambda$. 

\begin{Defn}[Residue classes and linkage classes]\leavevmode\newline
\begin{enumerate}
\item\vskip-15pt
Suppose that $\s,\t\in T(\Lambda)$. Then $\s$ and $\t$ are in the
same \textsf{residue class}, and we write $\s\approx\t$,
if $r_\s(i)=r_\t(i)$, for
$1\le i\le M$.
\item Suppose that $\lambda,\mu\in\Lambda$. Then $\lambda$ and $\mu$
are \textsf{residually linked}, and we write $\lambda\sim\mu$, if
there exist elements $\lambda_0=\lambda,\lambda_1,\dots,\lambda_r=\mu$
and elements $\s_j,\t_j\in T(\lambda_j)$ such that
$\s_{j-1}\approx\t_j$, for $i=1,\dots,r$.
\end{enumerate}
\end{Defn}

It is easy to see that~$\approx$ is an equivalence relation on
$T(\Lambda)$ and that~$\sim$ is an equivalence relation on~$\Lambda$.
If $\s\in T(\Lambda)$ let $\T_\s\in T(\Lambda)/\approx$ be its residue
class. If $\T$ is a residue class let $\T(\lambda)=\T\cap T(\lambda)$,
for $\lambda\in\Lambda$. By (\ref{L_K}), the residue classes
$T(\Lambda)/\approx$ parameterize the irreducible $\L_k$--modules.

Let $\T$ be a residue class $T(\Lambda)$ and define
$$F_\T = \sum_{\t\in\T}F_\t.$$ 
By definition, $F_\T$ is an element of $A_K$. We claim that, in fact,
$F_\T\in A_R$.

The following argument is an adaptation of Murphy's proof of
\citeAM[Theorem~2.1]{M:Nak}.

\begin{Lemma}\label{reduction}
Suppose that $\T$ is a residue equivalence class in $T(\Lambda)$. Then 
$F_\T$ is an idempotent in $A_R$.
\end{Lemma}

\begin{proof} We first note that $F_\T$ is an idempotent in~$A_K$
because it is a linear combination of orthogonal idempotents by
Theorem~\ref{idempotents}(a). The hard part is proving that 
$F_\T\in A_R$.

Fix an element $\t\in\T(\mu)$, where $\mu\in\Lambda$, and define
$$F_\t'=\prod_{i=1}^M\prod_{\substack{c\in\CC\\\bar c\ne r_\t(i)}}
              \frac{L_i-c}{c_\t(i)-c}.$$
Then $F_\t'\in A_R$ since, by assumption, $c_\t(i)-c$ is 
invertable in $R$ whenever $r_\t(i)\ne\bar c$. Observe that the 
numerator of $F_\t'$ depends only on $\T$ whereas the denominator
$d_\t=\prod_{i=1}^M\prod_{\bar c\ne r_\t(i)} (c_\t(i)-c)$ of~$F_\t'$ 
depends on~$\t$.
Let $\s\in T(\lambda)$. Then, by Proposition~\ref{vanishing}(d) and
Theorem~\ref{idempotents}(a),
$$F_\t' F_\s=\begin{cases}
         \frac{d_\s}{d_\t}F_\s,&\If\s\in\T,\\
         0,&\Otherwise.
\end{cases}$$ 
Consequently,
$F_\t'=\sum_{\lambda\in\Lambda}
   \sum_{\s\in \T(\lambda)} \frac{d_\s}{d_\t}F_\s$, 
by Theorem~\ref{idempotents}(c). 

Now, if $\s\in\T(\lambda)$ then $d_\s\equiv d_\t\pmod\pi$ since
$\s\approx\t$. Therefore, $1-\frac{d_\s}{d_\t}$ is a
non--zero element of $\pi$ since $d_\s\ne d_\t$ (as the
JM--elements separate $T(\Lambda)$ over $R$). Let $e_\s\in R$ be the
denominator of $F_\s$ and choose $N$ such that $e_\s\in\pi^N$, for all
$\s\in\T$. Then
$\big(1-\frac{d_\s}{d_\t}\big)^N\frac1{e_\s}\in R$, so that
$\big(1-\frac{d_\s}{d_\t}\big)^NF_\s\in A_R$, for all $\s\in\T$.
We now compute
\begin{align*}
\big(F_\T-F_\t'\big)^N
  &=\(\sum_{\lambda\in\Lambda}\sum_{\s\in\T(\lambda)}
        \big(1-\frac{d_\s}{d_\t}\big)F_\s\)^N\\
  &=\sum_{\lambda\in\Lambda}\sum_{\s\in\T(\lambda)}
        \big(1-\frac{d_\s}{d_\t}\big)^NF_\s,
\end{align*}
where the last line follows because the $F_\s$ are pairwise orthogonal
idempotents in~$A_K$. Therefore, $(F_\T-F_\t')^N\in A_R$. 

To complete the proof we evaluate $(F_\T-F_\t')^N$ directly. First, by 
Theorem~\ref{idempotents}(a), 
$$F_\t' F_\T=\sum_{\lambda\in\Lambda}\sum_{\s\in\T(\lambda)}
         \frac{d_\s}{d_\t}F_\s F_\T
      =\sum_{\lambda\in\Lambda}\sum_{\s\in\T(\lambda)}\frac{d_\s}{d_\t}F_\s
      =F_\t'.$$
Similarly, $F_\T F_\t'=F_\t'$. Hence, using the binomial theorem, we
have
\begin{align*}
(F_\T-F_\t')^N&=\sum_{i=0}^N(-1)^i\tbinom Ni(F_\t')^iF_\T^{N-i}\\
       &=F_\T+\sum_{i=1}^N(-1)^i\tbinom Ni (F_\t')^i\\
       &=F_\T+(1-F_\t')^N-1.
\end{align*}
Hence, $F_\T=(F_\T-F_\t')^N-(1-F_\t')^N+1\in A_R$, as required.
\end{proof}

By the Lemma, $F_\T\in A_R$. Therefore, we can reduce $F_\T$ modulo
$\pi$ to obtain an element of $A_k$. Let $G_\T=\bar{F_\T}\in A_k$ be
the reduction of $F_\T$ modulo~$\pi$. Then $G_\T$ is an idempotent
in~$A_k$.

Recall that if $\s\in T(\Lambda)$ then $\T_s$ is its residue class.

\begin{Defn} Let $\T$ be a residue class of $T(\Lambda)$. 
\begin{enumerate}
\item Suppose that $\s,\t\in\T(\lambda)$. Define 
$g^\lambda_{\s\t}=G_{\T_\s} a^\lambda_{\s\t}G_{\T_\t}\in A_k$.
\item Suppose that $\Gamma\in\Lambda/\sim$ is a residue linkage class in
$\Lambda$. Let $A_k^\Gamma$ be the subspace of~$A_k$ spanned by
$\set{g^\lambda_{\s\t}|\s,\t\in T(\lambda)\And\lambda\in\Gamma}$.
\end{enumerate}
\end{Defn}

Note that $G_\T^*=G_\T$ and that
$\big(g^\lambda_{\s\t}\big)^*=g^\lambda_{\t\s}$, for all 
$\s,\t\in T(\lambda)$ and $\lambda\in\Lambda$. By
Theorem~\ref{idempotents}, if $\mathbb S$ and $\T$ are residue classes
in $T(\Lambda)$ then $G_{\mathbb S}G_\T=\delta_{\mathbb S\T}G_\T$.

\begin{Prop}\label{G vanishing}
Suppose that $\s,\t\in T(\lambda)$, for some
$\lambda\in\Lambda$, that $\u\in T(\Lambda)$ and fix~$i$ with 
$1\le i\le M$. Let $\T\in\T(\Lambda)/\approx$. Then, in~$A_k$,
\begin{multicols}{2}
\begin{enumerate}
    \item $L_ig^\lambda_{\s\t}=r_\s(i)g^\lambda_{\s\t}$,
    \item $g^\lambda_{\s\t}L_i=r_\t(i)g^\lambda_{\s\t}$,
    \item $G_\T g^\lambda_{\s\t}=\delta_{\T_\s\T}\,g^\lambda_{\s\t}$,
    \item $g^\lambda_{\s\t}G_\T=\delta_{\T\T_\t}\,g^\lambda_{\s\t}$.
\end{enumerate}
\end{multicols}
\end{Prop}

We can now generalize the seminormal basis of the previous section to
the algebra~$A_k$.

\begin{Theorem}\label{blocks}
Suppose that $A_R$ has a family of JM--elements which separate
$T(\Lambda)$ over~$R$.
\begin{enumerate}
\item $\set{g^\lambda_{\s\t}|\s,\t\in T(\lambda)\And\lambda\in\Lambda}$
is a cellular basis of $A_k$.
\item Let $\Gamma$ be a residue linkage class of $\Lambda$. Then
$A_k^\Gamma$ is a cellular algebra with cellular basis
$\set{g^\lambda_{\s\t}|\s,\t\in T(\lambda)\And\lambda\in\Gamma}$.
\item The residue linkage classes decompose $A_k$ into a direct sum of 
cellular subalgebras; that is,
$$A_k=\bigoplus_{\Gamma\in\Lambda/\sim} A_k^\Gamma.$$
\end{enumerate}
\end{Theorem}

\begin{proof}Let $\Gamma$ be a residue linkage class in $\Lambda$ and
suppose that $\lambda\in\Gamma$. Then, exactly as in the proof of
Lemma~\ref{f basis}(a), we see that if $\s,\t\in T(\lambda)$ then
$g^\lambda_{\s\t}=a^\lambda_{\s\t}$ plus a linear combination of more
dominant terms. Therefore, the elements $\{g^\lambda_{\s\t}\}$ are
linearly independent because $\{a^\lambda_{\s\t}\}$ is a basis
of~$A_k$. Hence, $\{g^\lambda_{\s\t}\}$ is a basis of $A_k$.  We
prove the remaining statements in the Theorem simultaneously.

Suppose that $\lambda,\mu\in\Lambda$ and that
$\s,\t\in\T(\lambda)$ and $\u,\v\in\T(\mu)$. Then
$$g^\lambda_{\s\t}g^\mu_{\u\v} 
    =G_{\T_\s}a^\lambda_{\s\t}G_{\T_\t}G_{\T_\u}a^\mu_{\u\v}G_{\T_\v}
    =\begin{cases}
          G_{\T_\s}a^\lambda_{\s\t}G_{\T_\t}a^\mu_{\u\v}G_{\T_\v},
                 &\If \t\approx\u\\
          0,&\Otherwise.
\end{cases}$$
Observe that $\t\approx\u$ only if $\lambda\sim\mu$. Suppose then
that $\lambda\sim\mu$ and let $\Gamma$ be the residue linkage class in
$\Lambda$ which contains $\lambda$ and $\mu$. Then, because
$\{a^\nu_{\w\x}\}$ is a cellular basis of $A_k$, we can write
\begin{align*}
a^\lambda_{\s\t}G_{\T_\t}a^\mu_{\u\v}
     &=\sum_{\substack{\nu\in\Lambda\\\nu\ge\lambda,\nu\ge\mu}}
       \sum_{\substack{\w,\x\in T(\nu)\\\w\gedom\s,\x\gedom\v}}
             r_{\w\x}g^\nu_{\w\x},
\intertext{for some $r_{\w\x}\in k$ such that if $\nu=\lambda$ then
$r_{\w\x}\ne0$ only if $\w=\s$, and if $\nu=\mu$ then $r_{\w\x}\ne0$
only if $\x=\v$. Therefore, using Proposition~\ref{G vanishing}, we have}
g^\lambda_{\s\t}g^\mu_{\u\v} 
     &=\sum_{\substack{\nu\in\Lambda\\\nu\ge\lambda,\nu\ge\mu}}
       \sum_{\substack{\w,\x\in T(\nu)\\\w\gedom\s,\x\gedom\v}}
             r_{\w\x}G_{\T_s}g^\nu_{\w\x}G_{\T_\v}\\
     &=\sum_{\substack{\nu\in\Gamma\\\nu\ge\lambda,\nu\ge\mu}}
       \sum_{\substack{\w,\x\in T(\nu)\\\w\gedom\s,\x\gedom\v}}
             r_{\w\x}g^\nu_{\w\x}.
\end{align*}
Consequently, we see that if $\lambda\sim\mu\in\Gamma$ then
$g^\lambda_{\s\t}g^\mu_{\u\v}\in A_k^\Gamma$; otherwise,
$g^\lambda_{\s\t}g^\mu_{\u\v}=0$. All of the statements in the Theorem
now follow.
\end{proof}

Arguing as in the proof of Theorem~\ref{idempotents}(a) it follows
that $G_\T=\sum r_{\s\t}g^\lambda_{\s\t}$, where $r_{\s\t}$ is
non--zero only if $\s,\t\in\T(\lambda)$ for some $\lambda\in\Lambda$.

We are not claiming in Theorem~\ref{blocks} that the subalgebras
$A_k^\Gamma$ of $A_k$ are indecomposable. We call the indecomposable
two--sided ideals of $A_k$ the \textsf{blocks} of~$A_k$. It is a
general fact that each irreducible module of an algebra is a
composition factor of a unique block, so the residue linkage classes
induce a partition of the set of irreducible $A_k$--modules. By the
general theory of cellular algebras, all of the composition factors of
a cell module are contained in the same block; see~\citeAM[3.9.8]{GL}
or \citeAM[Cor.~2.22]{M:Ulect}. Hence, we have the following.

\begin{Cor}\label{necessary}
Suppose that $A_R$ has a family of JM--elements which separate
$T(\Lambda)$ over~$R$ and that $\lambda,\mu\in\Lambda$. Then
$C(\lambda)$ and $C(\mu)$ are in the same block of~$A_k$ only 
if~$\lambda\sim\mu$.
\end{Cor}

Let $\Gamma\in\Lambda/\sim$ be a residue linkage class. Then
$\sum_{\lambda\in\Gamma}F_\lambda\in A_R$ by Lemma~\ref{reduction} and
Theorem~\ref{idempotents}(b). Set 
$G_\Gamma=\bar{\sum_{\lambda\in\Gamma}F_\lambda}\in A_k$. The following
result is now immediate from Theorem~\ref{blocks} and
Theorem~\ref{idempotents}.

\begin{Cor}\label{G idempotents}
Suppose that $A_R$ has a family of JM--elements which
separate $T(\Lambda)$ over~$R$. 
\begin{enumerate}
\item Let $\Gamma$ be a residue linkage class. Then $G_\Gamma$ is a
central idempotent in $A_k$ and the identity element of the subalgebra
$A_k^\Gamma$. Moreover, 
$$A_k^\Gamma=G_\Gamma A_k G_\Gamma\cong\End_{A_k}(A_k G_\Gamma).$$
\item $\set{G_\Gamma|\Gamma\in\Lambda/\sim}$ and 
$\set{G_\T|\T\in T(\Lambda)/\approx}$ are complete sets of pairwise
orthogonal idempotents of~$A_k$. In particular,
$$1_{A_k} = \sum_{\Gamma\in\Lambda/\sim}G_\Gamma 
          = \sum_{\T\in T(\Lambda)/\approx} G_\T.$$
\end{enumerate}
\end{Cor}

Observe that the right ideals $G_\T A_k$ are projective
$A_k$--modules, for all $\T\in T(\Lambda)/\approx$. Of course, these
modules need not (and, in general, will not) be indecomposable.

Let $\RR(i)=\set{\bar c|c\in\CC(i)}$, for $1\le i\le M$. If $\T$ is a
residue class in $T(\Lambda)$ then we set $r_\T(i)=r_\t(i)$,
for $\t\in\T$ and $1\le i\le M$.

\begin{Cor}\label{G minimum poly}
Suppose that $A_R$ has a family of JM--elements which
separate $T(\Lambda)$ over~$R$. Then
$$L_i=\sum_{\T\in T(\Lambda)/\approx}r_\T(i)G_\T$$
and $\prod_{r\in\RR(i)}(L_i-r)$ is the minimum polynomial for $L_i$
acting on $A_k$.
\end{Cor}

\begin{proof} That
$L_i=\sum_{\T\in T(\Lambda)/\approx}r_\T(i)G_\T$
follows from Corollary~\ref{G idempotents}(b) and 
Proposition~\ref{G vanishing}. For the second claim,
for any $\s,\t\in T(\lambda)$ we have that
$$\prod_{r\in\RR(i)}(L_i-r)\cdot g^\lambda_{\s\t}=0$$
by Proposition~\ref{G vanishing}, so that
$\prod_{r\in\RR(i)}(L_i-r)=0$ in $A_k$. If we omit a factor
$(L_i-r_0)$ from this product then
$\prod_{r\ne r_0}(L_i-r)g^\lambda_{\s\t}\ne0$ whenever 
$\s,\t\in T(\lambda)$ and $r_0=r_\s(i)$. Hence, the product over
$\RR(i)$ is the minimum polynomial of $L_i$.
\end{proof}

As our final general result we note that the new cellular basis of
$A_k$ gives us a new `not quite orthogonal' basis for the cell modules
of $A_k$. Given $\lambda\in\Lambda$ fix $\s\in T(\lambda)$
and define $g^\lambda_\t=g^\lambda_{\s\t}+A_k^\lambda$ for $\t\in
T(\lambda)$.

\begin{Prop}
Suppose that $A_R$ has a family of JM--elements which
separate $T(\Lambda)$ over~$R$. Then 
$\set{g^\lambda_\t|\t\in T(\lambda)}$ is a basis of $C(\lambda)$.
Moreover, if $\t,\u\in T(\lambda)$ then
$$\<g^\lambda_\t,g^\lambda_\u\>_\lambda
       =\begin{cases}
	    \<a^\lambda_\t,g^\lambda_\u\>_\lambda,&\If\t\approx\u,\\
	    0,&\If\t\not\approx\u.
	\end{cases}$$
\end{Prop}

\begin{proof}
    That $\set{g^\lambda_\t|\t\in T(\lambda)}$ is a basis of
    $C(\lambda)$ follows from Theorem~\ref{blocks} and the argument of
    Lemma~\ref{f basis}(a). For the second
    claim, if $\t,\u\in T(\lambda)$ then
    $$\<g^\lambda_\t,g^\lambda_\u\>_\lambda 
            =\<a^\lambda_\t G_{\T_\t},g^\lambda_\u\>_\lambda
            =\<a^\lambda_\t,g^\lambda_\u G_{\T_\t}\>_\lambda$$
    by the associativity of the inner product since
    $G_{\T_\t}^*=G_{\T_\t}$. The result now
    follows from Proposition~\ref{G vanishing}(d).
\end{proof}

In the semisimple case Theorem~\ref{Gram det} reduces the Gram
determinant of a cell module to diagonal form. This result reduces it
to block diagonal form. Murphy has considered this block
decomposition of the Gram determinant for the Hecke algebras of
type~$A$~\citeAM{Murphy:ModGram}.

We now apply the results of this section to give a basis for the
blocks of several of the algebras considered in section~2.

\begin{Theorem}\label{application}
Let $k$ be a field and suppose that $A_R$ is one of the following 
algebras:
\begin{enumerate}
\item the group algebra $R\Sym_n$ of the symmetric group;
\item the Hecke algebra $\H_{R,q}(\Sym_n)$ of type $A$;
\item the Ariki--Koike algebra $\H_{R,q,\mathbf u}$ with $q\ne1$;
\item the degenerate Ariki--Koike algebra ${\mathcal H}_{R,\mathbf v}$;
\end{enumerate}
Then $A$ has a family of JM--elements which separate
$T(\Lambda)$ over~$R$ and Theorem~\ref{blocks} gives a basis for the
block decomposition of~$A_k$ into a direct sum of indecomposable
subalgebras.
\end{Theorem}

The cellular bases and the families of JM--elements for each of these
algebras are given in the examples of Section~2. As
$k\Sym_n\cong\H_{k,1}(\Sym_n)$, we use the Murphy basis for the
symmetric group. Note that the Hecke algebras of type $A$ should not
be considered as the special case $r=1$ of the Ariki--Koike algebras
because the JM--elements that we use for these two algebras are
different. Significantly, for the Ariki--Koike case we must assume
that $q\ne1$ as the JM--elements that we use do not separate
$T(\Lambda)$ over~$R$ when~$q=1$.

Before we can begin proving this result we need to describe how to
choose a modular system $(R,K,k)$ for each of the algebras above. In
all cases we start with a field $k$ and a non--zero element $q\in k$
and we let $R$ be the localization of the Laurent polynomial ring
$k[t,t^{-1}]$ at the maximal ideal generated by $(q-t)$. Then~$R$ is
discrete valuation ring with maximal ideal
$\pi$ generated by the image of~$(q-t)$ in~$R$. By construction,
$k\cong R/\pi$ and $t$ is sent to $q$ by the natural map
$R\longrightarrow k=R/\pi$. Let $K$ be the field of fractions of~$R$. 

First consider the case of the Hecke algebra $\H_{k,q}(\Sym_n)$.
As we have said, this includes the symmetric group as the
special case $q=1$. We take $A_R=\H_{R,t}(\Sym_n)$,
$A_K=\H_{K,t}(\Sym_n)$, and 
$A_k=\H_{R,t}(\Sym_n)\otimes_R k$. Then $\H_{K,t}(\Sym_n)$ is
semisimple and $\H_{k,q}(\Sym_n)\cong\H_{R,t}(\Sym_n)\otimes_R k$. 

Next, consider the Ariki--Koike algebra $\H_{k,q,\mathbf u}$ with
parameters $q\ne0,1$ and $\mathbf u=(u_1,\dots,u_m)\in k^m$. Let
$v_s=u_s+(q-t)^{ns}$, for $s=1,\dots,m$, and set $\mathbf
v=(v_1,\dots,v_m)$. We consider the triple of algebras
$A_R=\H_{R,t,\mathbf v}$, $A_K=\H_{K,t,\mathbf v}$ and
$A_k=\H_{k,q,\mathbf u}$. Once again, $A_K$ is semisimple and
$A_k\cong A_R\otimes_R k$. The case of the degenerate Ariki--Koike
algebras is similar and we leave the details to the reader.

The indexing set $\Lambda$ for each of the algebras considered in
Theorem~\ref{application} is the set of $m$--multipartitions of $n$,
where we identify the set of $1$--multipartitions with the set of
partitions. If $\lambda$ is an $m$--multipartition let $[\lambda]$ be
the diagram of $\lambda$; that is,
$$[\lambda]=\set{(s,i,j)|1\le s\le r\And 1\le j\le\lambda^{(s)}_i}.$$
Given a \textit{node} $x=(s,i,j)\in[\lambda]$ we define its content to be
$$c(x)=\begin{cases}
 [j-i]_t,&\If A_R=\H_{R,t}(\Sym_n),\\
 v_st^{j-i},&\If A_R=\H_{R,t,\mathbf v},\\
 v_s+(j-i),&\If A_R={\mathcal H}_{R,\mathbf v}.
\end{cases}$$
We set $\CC_\lambda=\set{c(x)|x\in[\lambda]}$ 
and $\RR_\lambda=\set{\bar{c(x)}|x\in[\lambda]}$. 

Unravelling the definitions, it is easy to see, for each of the
algebras that we are considering, that  if $\lambda\in\Lambda$ and
$\t\in T(\lambda)$ then $\CC_\lambda=\set{c_\t(i)|1\le i\le M}$.

To prove Theorem~\ref{application} we need to show that the residue
linkage classes correspond to the blocks of each of the algebras
above. Hence, Theorem~\ref{application} is a Corollary of the
following Proposition.

\begin{Prop}\label{residue classes}
    Let $A$ be one of the algebras considered in
    Theorem~\ref{application}. Suppose that $\lambda,\mu\in\Lambda$.
    The following are equivalent:
    \begin{enumerate}
	\item $C(\lambda)$ and $C(\mu)$ belong to the same block
	    of~$A_k$;
	\item $\lambda\sim\mu$;
	\item $\RR_\lambda=\RR_\mu$.
    \end{enumerate}
\end{Prop}

\begin{proof}
First suppose that $C(\lambda)$ and $C(\mu)$ are in the same block.
Then $\lambda\sim\mu$ by Corollary~\ref{necessary}, so that (a) implies
(b). Next, if (b) holds then, without loss of generality, there exist
$\s\in T(\lambda)$ and $\t\in T(\mu)$ with $\s\approx\t$; however,
then $\RR_\lambda=\RR_\mu$. So, (b) implies (c). The implication
`(c) implies (a)' is the most difficult, however, the blocks of all of
the algebras that we are considering have been classified and the
result can be stated uniformly by saying that the cell modules
$C(\lambda)$ and $C(\mu)$ belong to the same block if and only if
$\RR_\lambda=\RR_\mu$; see \citeAM{LM:blocks} for
$\H_{k,q}(\Sym_n)$ and the Ariki--Koike
algebras, and \citeAM{Brundan:degenCentre} for the degenerate
Ariki--Koike algebras. Therefore, (a) and (c) are equivalent. This
completes the proof.
\end{proof}

As a consequence we obtain the block decomposition of the Schur algebras.
Let $\Lambda_{m,n}$ be the set of $m$--multipartitions of $n$ and let
$S_{R,t,\mathbf v}(\Lambda_{m,n})$ be the corresponding cyclotomic
$q$--Schur algebra~\citeAM{DJM:cyc}, where $t$ and $\mathbf v$ are as
above.

\begin{Cor}
Let $k$ be a field and suppose that $A$ is one of the following 
$k$--algebras:
\begin{enumerate}
\item the $q$--Schur algebra $S_{R,q}(n)$;
\item the cyclotomic $q$--Schur $S_{R,t,\mathbf v}(\Lambda_{m,n})$
    algebra with $q\ne1$.
\end{enumerate}
Then $A$ has a family of JM--elements which separate $T(\Lambda)$
over~$R$ and Theorem~\ref{blocks} gives a basis for the block
decomposition of~$A_k$ into a direct sum of indecomposable
subalgebras.
\end{Cor}

\begin{proof}
Once again it is enough to show that two cell modules $C(\lambda)$ and
$C(\mu)$ belong to the same block if and only if $\lambda\sim\mu$.  By
Schur--Weyl duality, the blocks of $S_{k,q}(n)$ are in bijection with
the blocks of $\H_{k,q}(n)$~\citeAM[5.37--5.38]{M:Ulect} and the
blocks of $S_{k,q,\mathbf u}(\Lambda_{m,n})$ are in bijection with the
blocks of $\H_{k,q,\mathbf u}$~\citeAM[Theorem~5.5]{m:cyclosurv}.
Hence the result follows from Proposition~\ref{residue classes}.
\end{proof}

It is well known for each algebra $A$ in Theorem~\ref{application} the
symmetric polynomials in the JM--element belong to the centre of~$A$.
As our final result we show that there is a uniform explanation of
this fact. If $A$ is an algebra we let $Z(A)$ be its centre.

\begin{Prop}\label{symmetric poly}
Suppose that $A$ has a family of JM--elements which separate
$T(\Lambda)$ over~$R$ and that for $\lambda\in\Lambda$ there exist
scalars $c_\lambda(i)$, for $1\le i\le M$, such that
$$\set{c_\lambda(i)|1\le i\le M}=\set{c_\t(i)|1\le i\le M},$$
for any $\t\in T(\lambda)$. Then any
symmetric polynomial in $L_1,\dots,L_M$ belongs to the centre 
of~$A_k$.
\end{Prop}

\begin{proof}
Suppose that $X_1,\dots,X_M$ are indeterminates over~$R$ and let
$p(X_1,\dots,X_M)\in R[X_1,\dots,X_M]$ be a symmetric polynomial.
Recall that $L_i=\sum_\t c_\t(i)F_\t$ in~$A_K$, by 
Corollary~\ref{minimum poly}.  Therefore, 
$$p(L_1,\dots,L_M)
       =\sum_{\t\in T(\Lambda)}p\big(c_\t(1),\dots,c_\t(M)\big)F_\t
       =\sum_{\lambda\in\Lambda}
       p\big(c_\lambda(1),\dots,c_\lambda(M)\big)F_\lambda.$$
The first equality follows because the $F_\t$ are pairwise orthogonal
idempotents by Theorem~\ref{idempotents}.  By
Theorem~\ref{idempotents}(c) the centre of $A_K$ is spanned by the
elements $\set{F_\lambda|\lambda\in\Lambda}$, so this shows that
$p(L_1,\dots,L_M)$ belongs to the centre of $A_K$.  However,
$p(L_1,\dots,L_M)$ belongs to $A_R$ so, in fact, $p(L_1,\dots,L_M)$
belongs to the centre of $A_R$. Now, $\bar{Z(A_R)}$ is contained in
the centre of $A_k$ and any symmetric polynomial over~$k$ can be
lifted to a symmetric polynomial over~$R$. Thus, it follows that the
symmetric polynomials in the JM--elements of~$A_k$ are central
in~$A_k$.
\end{proof}

All of the algebras in Theorem~\ref{application} satisfy the
conditions of the Proposition because, using the notation above, if
$\t\in T(\lambda)$ then $\CC_\lambda=\set{c_\t(i)|1\le i\le M}$ for
any of these algebras. Notice, however, that the (cyclotomic) Schur
algebras considered in section~2 and the Brauer and BMW algebras do
not satisfy the assumptions of Proposition~\ref{symmetric poly}.


\bigskip
\centerline{\textit{Acknowledgements}}
I thank Marcos Soriano for many discussions about seminormal forms of
Hecke algebras and for his detailed comments and suggestions on this
paper.  This paper also owes a debt to Gene Murphy as he pioneered the
use of the Jucys--Murphy elements in the representation theory of
symmetric groups and Hecke algebras.


\newpage
\setcounter{section}{0}
\thispagestyle{firstpage}
\makeatletter
\title{Appendix. Constructing idempotents from triangular actions}
\begin{center}\baselineskip14pt\relax
{\bfseries\uppercasenonmath\@title\@title}\\
\bigskip
\textsc{Marcos Soriano}\footnote{The author thanks Andrew Mathas for his 
generosity and hospitality.}\\
\bigskip
\textit{Fachbereich Mathematik, Im Welfengarten 1,\\
        Universit\"at Hannover, Deutschland}\\
E-mail address: \texttt{soriano@math.uni-hannover.de}
\end{center}
\markboth{\uppercase{Marcos Soriano}}%
         {\uppercase{Constructing idempotents from triangular actions}}

\begin{abstract}
We give a general construction of a complete set of orthogonal
idempotents starting from a set of elements acting in an (upper)
triangular fashion. 
The construction is inspired in the Jucys--Murphy
elements (in their various appearances in several cellular algebras).  
\end{abstract}

\@setabstract
\makeatother

\section{Triangular actions: setup and notation.} \label{sec:1}
The construction of idempotents presented here is based only on matrix
arithmetic. 
However, whenever possible, we will mention the more suggestive 
notation from combinatorial representation theory.

Let $\Lambda$ be an $R$-algebra, where $R$ is an arbitrary integral
domain. 
The starting point is a representation $\rho$ of $\Lambda$ via matrices
over $R$, that is, an $R$-free (left) $\Lambda$-module $M$. 
Let $d$ be the $R$-rank of $M$ and set $\thi{d} := \{1,\ldots,d\}$. 

\begin{Remark} 
Until section~\ref{sec:5} we will not make any additional assumptions on 
$R$ or $\Lambda$. 
We have in mind such examples as $\Lambda$ being a cellular $R$-algebra and 
$M$ a single cell (``Specht'') module $M$, which would give rise to 
``Young's Orthogonal Form'' for $M$, as well as the case $M = \Lambda$ itself, 
e.g.\ for questions of semisimplicity.
\end{Remark}

Assume that with respect to a certain basis (of ``tableaux'') 
\[ \tcal := \{t_1, \ldots, t_d\}\subset M \]
there is a finite set of elements 
$\lcal := \{L_1, \ldots, L_n\}\subset\Lambda$ (the
``Jucys--Murphy'' elements) acting in an upper triangular way, that is,
\[ \rho(L_i) = \begin{pmatrix}
r_i^1&*&\cdots&*\\
0&r_i^2&\ddots&\vdots\\
\vdots&\ddots&\ddots&*\\
0&\cdots&0&r_i^d
\end{pmatrix}\,, \quad\forall i\in\thi{n}
\]
for certain diagonal entries $\{r_i^j\}$, $i\in\thi{n}, j\in\thi{d}$
(the ``residues'' or ``contents''). 
Call 
\[ (r_1^j, r_2^j, \ldots, r_n^j) \]
the \emph{residue sequence} corresponding to the basis element $t_j$. 
>From now on, we identify $L_i$ with its representing matrix,
thus suppressing $\rho$.
Note that we do \emph{not} make any assumption on $\<\lcal\>$ 
being central in $\Lambda$ or that $\lcal$ consists of pairwise commuting 
elements.

Finally, we need some notation for matrices. 
We denote by $\{E_{ij}\}_{i,j\in\underline{\textbf{d}}}$ the canonical matrix 
units basis of $\Mat_{d}(R)$, whose elements multiply according to 
$E_{ij}E_{kl} = \delta_{jk}E_{il}$.
The subring of $\Mat_{d}(R)$ consisting of upper triangular matrices contains
a nilpotent ideal with $R$-basis $\{E_{ij}\}_{1\le i<j\le d}$ which we denote
by $\ncal$.
We define the support of a matrix $A = (a_{ij})\in\Mat_{d}(R)$ in the
obvious way,
\[ \supp(A) := \{ (i,j)\in\thi{d}\times\thi{d}\,|\,a_{ij}\not= 0 \}\,. \]
To any $i\in\thi{d}$ we associate the following subset of $\thi{d}^{2}$:
\[ \ufrk_{i} := \{ (k,l)\in\thi{d}^{2}\,|\, k\le i\le l \}, \]
and extend this definition to any non--empty subset $J\subseteq\thi{d}$ via
$\ufrk_{J} := \bigcup_{i\in J}\ufrk_{i}$.
If $J$ is non--empty then a matrix $A$ has \emph{shape} $J$ if 
$\supp(A)\subseteq\ufrk_{J}$ and the sequence 
$(a_{ii})_{i\in\underline{\textbf{d}}}$ of diagonal entries is the 
characteristic function of the subset $J$, that is,
\[ a_{ii} = \left\{\begin{array}{ll}
1\,,&\text{if }i\in J\\
0\,,&\text{if }i\notin J\,.\\
\end{array}\right. \]
In particular, $A\in\sum_{i\in J}E_{ii}+\ncal$ and $A$ is upper triangular. 
For example, the matrices of shape $\{i\}$ have the form
\[\left(\begin{smallmatrix}
0&\cdots&0&*&*&\cdots&*\\
&\ddots&\vdots&\vdots&&&\vdots\\
&&0&*&&&*\\
&&&1&*&\cdots&*\\
&&&&0&\cdots&0\\
&&&&&\ddots&\vdots\\
&&&&&&0
\end{smallmatrix}\right)\,. \] 
\section{Cayley--Hamilton and lifting idempotents.}  \label{sec:2}
Let us pause to consider a single upper triangular matrix
\[ Z = \begin{pmatrix}
\zeta_{1}&*&\cdots&*\\
0&\zeta_{2}&\ddots&\vdots\\
\vdots&\ddots&\ddots&*\\
0&\cdots&0&\zeta_{d}
\end{pmatrix}\in\Mat_{d}(R)\,. \]
Note that by the Cayley--Hamilton theorem, the matrix $Z$ satisfies the
polynomial $\prod_{i=1}^{d}(X-\zeta_{i})$.
Assume that $Z$ has shape $J$ for some non--empty 
$J\subseteq\thi{d}$ of cardinality $k=|J|$.
Then $Z$ satisfies the polynomial $(X-1)^{k}\cdot X^{d-k}$. 
What if $k = 1$? 
Then the Cayley--Hamilton equation for $Z$ reads
\[ 0=Z^{d-1}\cdot(Z-1) \:\Leftrightarrow\: Z^{d} = Z^{d-1}\,. \]
This implies (by induction) $Z^{d+j} = Z^{d}$ for all $j\ge 1$.
In particular, the element $F := Z^{d}$ is an idempotent. 

Of course, this is just a special case of ``lifting'' idempotents, and
can be extended (cf.\ \citeMS{Feit}, Section I.12) 
to the following ring theoretical version (Lemma \ref{l1}). 
We introduce some notation first. \\
Let $N\ge 2$ be a natural number (corresponding to the nilpotency degree in
Lemma \ref{l1}; for $N=1$ there is nothing to do).
Consider the following polynomial in two (commuting) indeterminates
\begin{eqnarray*}
(X+Y)^{2N-1} & = & \sum_{i=0}^{2N-1}\binom{2N-1}{i}X^{i}Y^{2N-1-i} \\
& = &
\sum_{i=0}^{N-1}\binom{2N-1}{i}X^{2N-1-i}Y^{i} +
\sum_{i=0}^{N-1}\binom{2N-1}{i}X^{i}Y^{2N-1-i} \\
& =: & \varepsilon_{N}(X,Y)+\varepsilon_{N}(Y,X) 
\end{eqnarray*}
(using the symmetry of the binomial coefficients). 
Note that $\varepsilon_{N}(X,Y)$ has \emph{integer} coefficients.    
Since $N > 1$,
\begin{equation}\label{eq1} 
\varepsilon_{N}(X,Y) \equiv X^{2N-1} \mod(XY)
\end{equation}
and $\varepsilon_{N}(X,Y) \equiv 0 \mod (X^{N})$. This implies that
\begin{equation}\label{eq2} 
\varepsilon_{N}(X,Y)\cdot\varepsilon_{N}(Y,X) \equiv 0 \mod (XY)^{N}\,.
\end{equation}
Specialise to a single indeterminate by setting $\varepsilon_{N}(X) :=
\varepsilon_{N}(X,1-X)$ and observe that
\begin{equation}\label{eq3} 1 = 1^{2N-1} = (X+(1-X))^{2N-1} =
\varepsilon_{N}(X)+\varepsilon_{N}(1-X)\,. 
\end{equation}
Now we are ready to formulate the
\begin{Lemma}\label{l1}
Let $\hcal$ be a ring and $\ncal$ a nilpotent two--sided ideal of 
nilpotency degree $N$. \\
If $x^{2} \equiv x \mod \ncal$, then $e := \varepsilon_{N}(x)$ is an 
idempotent with $e \equiv x \mod \ncal$.  
\end{Lemma}
\begin{proof} Note that $x^{2} \equiv x \mod \ncal \Leftrightarrow
x-x^{2}=x(1-x)\in \ncal$, implying
\[ e = \varepsilon_{N}(x) \equiv x^{2N-1} \equiv x \mod \ncal \]
by equation (\ref{eq1}). 
On the other hand, combining equations (\ref{eq2}) and (\ref{eq3}),
\[ e-e^{2} = e(1-e) = \varepsilon_{N}(x)\varepsilon_{N}(1-x) \equiv 0 \mod
(x(1-x))^{N}\,. \]
But $(x(1-x))^{N}\in \ncal^{N}=0$, thus the equality $e-e^{2}=0$ holds and $e$
is indeed an idempotent.
\end{proof}
\section{The separating condition and directedness.}  \label{sec:3}
We consider first a simple version of the idempotent construction
that is relevant to semisimplicity questions.
For $i\in\thi{d}$ we denote by $\mis{i}$ the set $\thi{d}\setminus\{i\}$. 
Let us assume now that \emph{for all} $i\in\thi{d}$ the following 
\emph{separating} condition is satisfied:
\[ (\scal)\qquad \forall j\in\mis{i}\quad \exists k\in\thi{n}
\text{ such that }r_{k}^{i}-r_{k}^{j}\in R^{\times}\,. \]
In particular, both residues are different.  
Of course, $k=k(j)=k(j,i)$ may not be unique, but we assume a 
fixed choice made for all possible pairs of indices. 
Then we define ($\One$ denotes the identity matrix)
\[ Z_{i} := \prod_{j\in\hat{{\rm\underline{\textbf{i}}}}}
\frac{L_{k}-r_{k}^{j}\One}{r_{k}^{i}-r_{k}^{j}}\,.\]
The product can be taken in \emph{any} order, the essential fact being
only to achieve that the matrix $Z_{i}$ is of the form $Z_{i} = E_{ii}+N_{i}$
for some upper triangular nilpotent matrix $N_{i}$.
Just note that for the $j$-th factor $F$ in the definition of $Z_{i}$ we have
\[ F_{ii} = \frac{r_{k}^{i}-r_{k}^{j}}{r_{k}^{i}-r_{k}^{j}} = 1 
\quad\text{ and }\quad 
F_{jj} = \frac{r_{k}^{j}-r_{k}^{j}}{r_{k}^{i}-r_{k}^{j}} = 0\,. \]
Now, using the observation of \S 2, we obtain a set of \emph{idempotents} 
$\ecal_{i} := Z_{i}^{d}$.
Our first assertion is
\begin{Lemma}\label{l2}
The idempotent $\ecal_{i}$ has shape $\{i\}$.
\end{Lemma}
\begin{proof}
Any matrix of the form $(E_{ii}+N)^{d}$ with $N\in\ncal$ has shape $\{i\}$.
To see this, use the non--commutative binomial expansion for 
$U = (E_{ii}+N)^d$, that is, express $U$ as a sum of terms 
$X_1\cdots X_d$, where $X_j\in\{E_{ii},N\}$.
If all $X_j = N$, we have the (only) summand of the form $N^d = 0$ 
(by nilpotency), with no contribution. 
Similarly, if all $X_j = E_{ii}$, we obtain one summand $E_{ii}$. \\
In the case when $X_1$ or $X_d$ equals~$E_{ii}$, and at least one factor
equals $N$, this summand has the 
appropriate form,
\[ \text{either }\quad
\left(\begin{smallmatrix}
0&\cdots&0&0&0&\cdots&0\\
&\ddots&\vdots&\vdots&\vdots&&\vdots\\
&&0&0&0&\cdots&0\\
&&&0&*&\cdots&*\\
&&&&0&\cdots&0\\
&&&&&\ddots&\vdots\\
&&&&&&0
\end{smallmatrix}\right)
\quad\text{ or }\quad
\left(\begin{smallmatrix}
0&\cdots&0&*&0&\cdots&0\\
&\ddots&\vdots&\vdots&\vdots&&\vdots\\
&&0&*&0&\cdots&0\\
&&&0&0&\cdots&0\\
&&&&0&\cdots&0\\
&&&&&\ddots&\vdots\\
&&&&&&0
\end{smallmatrix}\right)\,.
\] 
Thus we are left with the summands having $X_1=X_d=N$ and $X_j=E_{ii}$ for 
some $1<j<d$. 
But the support of any matrix in $\ncal E_{ii}\ncal$ is contained in
the set $\{ (k,s)\in\thi{d}^{2}\,|\, k<i<s \}$, as one sees
by matrix unit gymnastics (running indices are underlined):
\begin{eqnarray*}
\big(\sum_{\ul{k}<\ul{j}}a_{kj}E_{kj}\big)\cdot E_{ii}\cdot 
\big( \sum_{\ul{r}<\ul{s}}b_{rs}E_{rs}\big) 
= \big(\sum_{\ul{k}<\ul{j}}a_{kj}E_{kj}\big)
\big(\sum_{i<\ul{s}}b_{is}E_{is}\big)
= \sum_{1\le\ul{k}<i<\ul{s}\le d}a_{ki}b_{is}E_{ks}\,.
\end{eqnarray*}
This finishes the proof of the lemma, as all summands add up to
give $U - E_{ii}$ nilpotent with support contained in $\ufrk_{i}$.
\end{proof}
Lemma \ref{l2} has an important consequence: the one--sided ``directed''
orthogonality of the obtained idempotents.
\begin{Defn}
Let $\hcal$ be an arbitrary ring. Call a finite set $\{e_{1},\ldots,e_{d}\}$
of idempotents in $\hcal$ \emph{directed}, if $e_{j}e_{i} = 0$ 
whenever $j > i$. 
\end{Defn}
\begin{Lemma}\label{l3}
The set of idempotents $\{\ecal_i\}_{i\in\underline{\textbf{d}}}$ 
is directed.
\end{Lemma}
\begin{proof}
Directedness is an immediate consequence of the fact that $\ecal_{i}$ 
has shape~$\{i\}$.
\end{proof}
\section{Gram--Schmidt orthogonalisation of directed systems of idempotents.}
\label{sec:4}
We can now proceed inductively and construct a complete set of orthogonal
idempotents out of $\{\ecal_{i}\}$. 
The inductive step goes as follows: 
\begin{Lemma}\label{l4}
Let $\hcal$ be an arbitrary ring.
Assume we are given two finite sets of idempotents in $\hcal$
(one of them possibly empty)
\[ \ebb = \{e_{1},\ldots, e_{k}\} \quad\text{ and }\quad
\fbb =\ \{f_{k+1}, f_{k+2}, \ldots, f_{d}\} \]
for some $k\ge 0$ with the following properties:
\begin{enumerate}
\item $\ebb$ consists of pairwise orthogonal idempotents\,,
\item $\fbb$ is directed\,,
\item $\ebb$ is orthogonal to $\fbb$, that is, $ef = 0 =fe$ for $e\in\ebb$,
$f\in\fbb$.
\end{enumerate}
Set $F := \sum_{i=1}^{k}e_{i}+f_{k+1}$. 
Then the sets of idempotents
\[ \widetilde{\ebb} = \{e_{1},\ldots, e_{k}, f_{k+1}\} \quad\text{ and }\quad
\widehat{\fbb} =\ \{(1-F)f_{k+2}, \ldots, (1-F)f_{d}\} \]
satisfy conditions $(a)$--$(c)$.
\end{Lemma}
\begin{proof}
First observe that $F$ is an idempotent, by orthogonality. 
If $j\ge k+2$ we have (by the orthogonality of~$\ebb$ and $\fbb$ and
the directedness of $\fbb$) that
\[ f_{j}\cdot F = f_{j}\cdot(e_{1}+\ldots+e_{k}+f_{k+1}) 
= \sum_{i=1}^{k}f_{j}e_{i} + f_{j}f_{k+1} = 0 + 0 = 0\,. \]
This implies that $\hat{f}_j := (1-F)f_{j}$ is an idempotent because
\begin{eqnarray*}
\hat{f}_{j}^{2} = (f_{j}-Ff_{j})(f_{j}-Ff_{j})
= f_{j}-Ff_{j}-\underbrace{f_{j}F}_{=0}f_{j}+F\underbrace{f_{j}F}_{=0}f_{j} 
= (1-F)f_{j} = \hat{f}_{j}.
\end{eqnarray*} 
Similarly, the set $\{\hat{f}_{s}\}_{k+2\le s\le d}$ 
is directed because for $j>i>k+1$
\begin{eqnarray*}
\hat{f}_{j}\cdot\hat{f}_{i} = (f_{j}-Ff_{j})(f_{i}-Ff_{i}) 
= \underbrace{f_{j}f_{i}}_{=0}-\underbrace{f_{j}F}_{=0}f_{i}-
F\underbrace{f_{j}f_{i}}_{=0}+F\underbrace{f_{j}F}_{=0}f_{i} = 0\,.
\end{eqnarray*}
Since $\ebb$ is orthogonal to $f_{k+1}$, $\widetilde{\ebb}$ consists
obviously of pairwise orthogonal idempotents. 
Thus, we are left with checking orthogonality between $f_{k+1}$ and 
$\widehat{\fbb}$.
Let $j\ge k+2$, then
\begin{eqnarray*}
f_{k+1}\cdot\hat{f}_{j} = f_{k+1}(1-F)f_{j} 
= \underbrace{f_{k+1}(1-f_{k+1})}_{=0}f_{j}
-\sum_{i=1}^{k}\underbrace{f_{k+1}e_{i}f_{j}}_{=0} = 0\,,
\end{eqnarray*}
as well as $\hat{f}_{j}\cdot f_{k+1} = (1-F)f_{j}f_{k+1} = 0$ by directedness.
\end{proof}

Thus, keeping the notations from \S 1 and \S 3, we obtain the following
\begin{Prop}\label{p1}
A set $\lcal = \{L_{1},\ldots,L_{n}\}$ of ``Jucys--Murphy operators''
satisfying the separating condition $(\scal)$ for all $i\in\thi{d}$ gives
rise to a complete set of orthogonal idempotents $\{e_{1},\ldots,e_{d}\}$. 
\end{Prop}
\begin{proof}
Starting from $\ebb = \emptyset$ and 
$\fbb = \{\ecal_{i} = Z_{i}^{d}\}_{i\in\underline{\textbf{d}}}$,
we obtain --- using Lemma \ref{l4} $d$ times ---  a set
$\{e_{1},\ldots,e_{d}\}$ of orthogonal idempotents. \\
Note that the idempotents $e_{i}$ have again shape $\{i\}$ 
(check this in the inductive step from Lemma \ref{l4} 
by considering the form of the matrix $\One-F$).
Completeness of the set $\{e_1,\ldots,e_d\}$ now follows easily, since we 
obviously have by Lemma \ref{l2}:
\[ e := e_{1}+\ldots+e_{d} = \One+N\,, \]
for some (upper triangular) nilpotent matrix $N$. 
Thus, $\One-e$ is an idempotent \emph{and} a nilpotent matrix, 
implying that $N=0$. 
\end{proof}
Note that the proof gives, at the same time, a practical 
\emph{algorithm} for \emph{constructing} the complete set of orthogonal
idempotents in question.
\section{Linkage classes.}  \label{sec:5}
>From now on, we assume that $R$ is a \emph{local} commutative ring with
maximal ideal $\mfrk$. 
This includes the case of $R$ being a field (when $\mfrk = 0$).

Fix $k\in\thi{n}$ and $j\in\thi{d}$. 
We may assume without loss of generality that not all residues 
$r_{k}^{i}$, $i\in\thi{d}$,
are zero (replace $L_{k}$ by $\One+L_{k}$ if necessary
\footnote{Note that this does not change the property of the considered set of 
Jucys--Murphy operators of being central in $\Lambda$ or, rather, 
consisting of pairwise commuting elements.}).
We say that $i\in\thi{d}$ is \emph{linked} to $j$ via $L_{k}$, 
if $r_{k}^{i}-r_{k}^{j}\in\mfrk$.
Set 
\[ \lbb_{k}(j) := \{ i\in\thi{d}\,|\, i\text{ is linked to }j\text{ via }
L_{k} \}\,. \]
Observe that $j\in\lbb_{k}(j)$ since $0\in\mfrk$.
\begin{Defn}
The \emph{linkage class} of $j\in\thi{d}$ with respect to 
$\lcal = \{L_{1},\ldots,L_{n}\}$ is the set
\[ \lbb(j) := \bigcap_{k\in\underline{\textbf{n}}}\lbb_{k}(j)\,. \]
\end{Defn}
\begin{Remark}
Linkage classes with respect to $\lcal$ partition the set $\thi{d}$
(of ``tableaux'') into, say, $l$ disjoint sets $J_{1},\ldots, J_{l}$.
In view of the fact that $R\setminus\mfrk = R^{\times}$, the assumption
of the separating condition $(\scal)$ from \S 3 for all $i\in\thi{d}$
just translates into the condition of all linkage classes being singletons. 
\end{Remark}
Consider a fixed linkage class $J$.
For all $j\in\thi{d}\setminus J$ we assume that a fixed choice of
$k\in\thi{n}$ and $i\in\thi{d}$ has been made such that
\[ r_{k}^{i}-r_{k}^{j}\in R^{\times} = R\setminus\mfrk\,. \]
Then we define
\[ Z_{J} := \prod_{j\notin J}
\frac{L_{k}-r_{k}^{j}\One}{r_{k}^{i}-r_{k}^{j}} \]
(the product can be taken in any order). Note that --- by 
Lemma \ref{l2} --- $Z_{J}^{d}$ has shape $J$.
\section{A general orthogonalisation algorithm for idempotents.}
\label{sec:6}
\begin{Prop}\label{p2}
A set $\lcal = \{L_{1},\ldots,L_{n}\}$ partitioning $\thi{d}$ into
$l$ linkage classes gives rise to a complete set $\{e_{1},\ldots,e_{l}\}$
of orthogonal idempotents. 
\end{Prop}
Again, the proof of the proposition provides an algorithm for constructing
$\{e_{1},\ldots,e_{l}\}$.
\begin{proof}
Let $J_{1},\ldots,J_{l}$ denote the linkage classes and set
$U_{i} := Z_{J_{i}}^{d}$, a matrix of shape $J_{i}$.
We start the orthogonalisation procedure by setting
$\ebb_{0} := \emptyset$ and 
$\fbb_{l} := \{\varepsilon_{d}(U_{i})\}_{1\le i\le l}$.
Note that $\fbb_{l}$ consists of idempotents by Lemma \ref{l1}.
Assuming that two sets of idempotents $\ebb_{k} = \{e_{1},\ldots,e_{k}\}$ 
(pairwise orthogonal) and $\fbb_{l-k} = \{f_{k+1},\ldots,f_{l}\}$ with
$\ebb_{k}$ orthogonal to $\fbb_{l-k}$ have been already constructed,
we set $\ebb_{k+1} := \ebb_{k}\cup\{f_{k+1}\}$ and have to modify
$\fbb_{l-k}$ appropriately.
The goal is that $\fbb_{l-k-1}$ consists of idempotents orthogonal to 
$\ebb_{k+1}$.

Set $F := \sum_{e\in\ebb_{k+1}}e$ and consider first
$\tilde{f}_{j} := \varepsilon_{d}\big((1-F)f_{j}\big)$ for all $j\ge k+2$.
Since $e(1-F) = 0$ for $e\in\ebb_{k+1}$ and 
$\varepsilon_{d}(X)\in X^{d}\cdot\zbb[X,Y]$, $\ebb_{k+1}$ is \emph{left}
orthogonal to the idempotent $\tilde{f}_{j}$, $j\ge k+2$.
Similarly, multiplication from the right by $(1-F)$ and application
of the polynomial $\varepsilon_{d}$ forces \emph{right} orthogonality
to hold, while keeping left orthogonality.
That is, the set $\fbb_{l-k-1} = \{\hat{f}_{k+2},\ldots,\hat{f}_{d}\}$ with
$\hat{f}_{j} := \varepsilon_{d}\big(\tilde{f}_{j}(1-F)\big)$ has
the desired properties.

Thus, after $l$ steps, we end up with an orthogonal set of idempotents
$\{e_{1},\ldots,e_{l}\}$. 
Observe that the inductive step described above does not change the
shape of the idempotents, implying that $e_{i}$ has shape
$J_{i}$ as the original idempotent $U_{i}$.
This fact, in addition to $J_{1},\ldots,J_{l}$ partitioning $\thi{d}$,
leads to the equation
\[ e_{1}+\ldots+e_{l} = \One-N \]
with $N$ a nilpotent and idempotent matrix, thus implying $N = 0$ and
the completeness of $\ebb_{l}$. 
\end{proof}
\begin{Remark}
Retracing all steps in the proof of Proposition \ref{p2}, we see
that the constructed idempotents $e_{i}$ belong to $R[L_{1},\ldots,L_{n}]$, 
the $R$-subalgebra of $\Lambda$ generated by $\lcal$.
Thus, if the elements from $\lcal$ do commute pairwise, this will still
hold for the set of idempotents $\ebb := \ebb_{l}$.

In particular, assuming that $\lcal$ is a set of \emph{central} Jucys--Murphy
elements for the module $M = \Lambda$, we obtain a set $\ebb$ of \emph{central}
orthogonal idempotents.
Thus, for example, the \emph{block} decomposition of $\Lambda$ in the
case of $R$ being a field must be a refinement of the decomposition
into linkage classes induced by $\lcal$.

We leave the adaptation of the presented methods to particular classes 
or examples for $\Lambda$, $R$, $M$ and $\lcal$ to the reader's needs.
\end{Remark}

\end{document}